\documentclass[a4paper,twoside]{amsart}
\usepackage[utf8]{inputenc}
\usepackage[T1]{fontenc}
\usepackage[british]{babel}
\usepackage{amssymb,amsmath,amsthm,amscd,amsfonts}
\usepackage{color}
\usepackage{verbatim}
\usepackage{lmodern}
\usepackage{enumerate}
\usepackage{hyperref}
\usepackage{tikz}
\usetikzlibrary{arrows.meta}

\theoremstyle{definition}
\newtheorem{thmA}{Theorem}

\newtheorem{corA}[thmA]{Corollary}
\newtheorem{thm}{Theorem}[section]
\newtheorem{defn}[thm]{Definition}
\newtheorem{lemma}[thm]{Lemma}
\newtheorem{prop}[thm]{Proposition}
\newtheorem{cor}[thm]{Corollary}

\newtheorem{rem}[thm]{Remark}

\newcommand\Z{\mathbf{Z}}
\newcommand\E{\mathbf{E}}
\newcommand\F{\mathbf{F}}
\newcommand\prob{\mathbf{P}}
\newcommand\C{\mathbf{C}}
\newcommand\R{\mathbf{R}}
\newcommand\N{\mathbf{N}}
\newcommand\x{\mathbf{x}}

\newcommand\SL{\mathrm{SL}}
\newcommand\Sp{\mathrm{Sp}}
\renewcommand\O{\mathrm{O}}

\newcommand\func[1]{\varphi(#1)}

\definecolor{darkgreen}{rgb}{0.0, 0.5, 0.0}

\author[de Laat]{Tim de Laat}
\thanks{T.dL.~was funded by the Deutsche Forschungsgemeinschaft -- Project-ID 427320536 -- SFB 1442, as well as under Germany’s Excellence Strategy EXC 2044 390685587, Mathematics M\"unster: Dynamics--Geometry--Structure}
\address{Department of Mathematics, Kiel University \newline Heinrich-Hecht-Platz 6, D-24118 Kiel, Germany}
\email{delaat@math.uni-kiel.de}

\author[de la Salle]{Mikael de la Salle}
\thanks{M.dlS.~was supported by the ANR project ANCG Project-ANR-19-CE40-0002}
\address{Université Lyon 1, Centrale Lyon, INSA Lyon, Université Jean Monnet, CNRS, ICJ UMR5208, 69622 Villeurbanne, France}
\email{delasalle@math.univ-lyon1.fr}

\title{Actions of higher rank groups on uniformly convex Banach spaces}

\date{\today}

\begin{document}

\begin{abstract}
    We prove that all isometric actions of higher rank simple Lie groups and their lattices on arbitrary uniformly convex Banach spaces have a fixed point. This vastly generalises a recent breakthrough of Oppenheim. Combined with earlier work of Lafforgue and of Liao on strong Banach property (T) for non-Archimedean higher rank simple groups, this confirms a long-standing conjecture of Bader, Furman, Gelander and Monod. As a consequence, we deduce that sequences of Cayley graphs of finite quotients of a higher rank lattice are super-expanders.
\end{abstract}

\maketitle

\section{Introduction}
Given a Banach space $E$, a topological group $G$ has property F$_E$ if every continuous action $G \curvearrowright E$ by affine isometries has a fixed point. This property was introduced by Bader, Furman, Gelander and Monod \cite{MR2316269} as a Banach space version of Kazhdan's property (T).

Probably the best-known examples of groups with property (T) are higher rank groups and lattices in such groups. Throughout this article, a higher rank group is a group of the form $G = \prod_{i=1}^n G_i (\F_i)$, where each $\F_i$ is an arbitrary local field and each $G_i$ is a Zariski connected (almost) $\F_i$-simple group with $\F_i$-rank $\geq 2$. If $\F_i=\R$, then $G_i(\F_i)$ is a connected simple Lie group with finite center and real rank $\geq 2$.

In \cite[Theorem B]{MR2316269}, Bader, Furman, Gelander and Monod proved that higher rank groups and their lattices satisfy property F$_E$ for all $L^p$-spaces $E$ and for a large class of subquotients $E$ of $L^p$-spaces, whenever $1 < p < \infty$. They conjectured a much stronger statement, which we prove to be true as a consequence of our main result (Theorem \ref{thm:higherrankgroupsFE}), which itself is a vast generalization of a recent breakthrough of Oppenheim, combined with earlier work of Lafforgue and of Liao. Before we state this result, recall that a Banach space $E$ is super-reflexive if it carries an equivalent uniformly convex norm.
\begin{thmA}{\cite[Conjecture 1.6.]{MR2316269}} \label{thm:BFGMconjecture}
    Higher rank groups and their lattices have property F$_E$ for every super-reflexive Banach space $E$.
\end{thmA}
To illustrate the strength of Theorem \ref{thm:BFGMconjecture}, recall that every locally compact, second countable group admits a proper (and hence fixed-point free) affine isometric action on an $L^{\infty}$-space, and even on a strictly convex, reflexive (but not super-reflexive) Banach space; see \cite{HaagerupPrzybyszewska2006}. On the other hand, it is known that every affine isometric action of a locally compact, second countable group with property (T) on an $L^1$-space has a fixed point \cite{MR2929085}.

The non-Archimedean case of \cite[Conjecture 1.6]{MR2316269} (i.e.~the case in which each $\F_i$ is a non-Archimedean local field) is a consequence of the groundbreaking work of Lafforgue \cite{MR2423763,MR2574023} and of Liao \cite{MR3190138} on strong Banach property (T). This powerful approach using strong property (T) has also lead to partial results in the Archimedean case of the conjecture \cite{MR3474958,MR3407190,MR3533271,MR2957217}, but the conjecture remained open.

Recently, Oppenheim established a breakthrough in the Archimedean case. Indeed, in \cite[Theorem 1.2]{Oppenheim2022}, he proved that the group $\SL_3(\Z)$ (and hence also the group $\SL_3(\R)$) has property (T$_E$) for all super-reflexive Banach spaces $E$. From this, it follows (see \cite[Theorem 1.5]{Oppenheim2022}) that for all $n \geq 4$, the group $\SL_n(\R)$ and its lattices have property F$_E$.

The aim of this article is to prove the following generalisation of this result.
\begin{thmA} \label{thm:higherrankgroupsFE}
    Let $G$ be a connected simple Lie group with real rank $\geq 2$ and finite center. Every continuous affine isometric action of $G$ and of any lattice in $G$ on a uniformly convex Banach space has a fixed point. In other words, $G$ and its lattices have property F$_E$ for every uniformly convex Banach space $E$.
\end{thmA}
From the point of view of the structure of $G$, the proof of this result follows the same strategy of proof as Oppenheim's result for $\SL_3(\Z)$. Concretely, it relies on analysis on several different (isomorphic) nilpotent subgroups of $G$. There are also important differences between Oppenheim's approach and our approach, the main one being that we directly consider isometric actions on Banach spaces rather than isometric representations. We give an outline of the strategy of proof of Theorem \ref{thm:higherrankgroupsFE} in Section \ref{subsec:strategy}, highlighting the similarities to and differences from Oppenheim's approach.

%We point out that, when considering the specific case of Hilbert spaces, our proof of Theorem \ref{thm:higherrankgroupsFE} is the first direct proof of property FH (the fixed point property for isometric actions on Hilbert spaces) that does not use any spectral theory.\footnote{I forget: is this really correct? Doesn't Oppenheim's proof also avoid any spectral theory?}

We explain how Theorem \ref{thm:higherrankgroupsFE} and the known non-Archimedean case (as recalled above) together imply Theorem \ref{thm:BFGMconjecture} in Section \ref{sec:conclusion}.

Theorem \ref{thm:BFGMconjecture} has interesting consequences in the direction of super-expanders and spectral gap, which follow by known methods that we recall in Section \ref{sec:superexpanders}. We also refer to that section for the precise definition of (super-)expander.

It is a deep open problem whether there exists an expander that is not a super-expander. Already constructing super-expanders has proved a difficult task, and only a few classes of examples are known, namely sequences of Cayley (or Schreier) graphs of finite quotients of groups with Banach property (T) \cite{MR2423763,MR3190138}, super-expanders constructed by means of the zigzag product \cite{MR3210176}, expanders constructed from a warped cone over an action of a group with strong Banach property (T) \cite{MR3956196,MR4245587,MR3912479}, and by the recent work of Oppenheim \cite{Oppenheim2022}, sequences of Cayley (or Schreier) graphs of finite quotients of any lattice in a simple Lie group locally containing $\SL_3(\R)$.

We highlight the following striking consequence of Theorem \ref{thm:BFGMconjecture}, which is a vast generalisation of the first and the fourth class of super-expanders mentioned above.
\begin{corA} \label{cor:superexpanders}
Let $\Gamma$ be a lattice in a higher rank group, let $S$ be a finite symmetric generating set of $\Gamma$, and let $(\Gamma_n)$ be a sequence of finite quotients of $\Gamma$ such that $|\Gamma_n| \to \infty$ for $n \to \infty$. Then the sequence $\mathrm{Cay}(\Gamma_n,S)$ of Cayley graphs is a super-expander.
\end{corA}

\subsection{Strategy of proof} \label{subsec:strategy}
We briefly outline the strategy of proof of Theorem \ref{thm:higherrankgroupsFE}. First, by well-known arguments that we recall in Section~\ref{sec:conclusion}, it suffices to prove Theorem~\ref{thm:higherrankgroupsFE} for actions of the rank $2$ simple Lie groups $\SL_3(\R)$ and $\Sp_4(\R)$, the main reason being that every connected simple higher rank Lie group has a closed subgroup locally isomorphic to one of these groups. This is done in Theorem~\ref{thm:FESL3} and Theorem~\ref{thm:FESp4}, respectively. The strategy in both cases is similar, but the details are significantly more involved for $\Sp_4(\R)$. Our proof is completely self-contained, but we use some of Oppenheim's ideas from \cite{Oppenheim2022}.

As mentioned above, Oppenheim proved the theorem for $\SL_n(\R)$ for $n \geq 4$. In fact, his main original result is not about affine actions of $\SL_{n \geq 4}(\R)$, but about linear actions of $\SL_3(\Z)$. Recall that a group has property (T$_E$) if every isometric representation on $E$ has spectral gap \cite[Definition 1.1]{MR2316269}; see Section~\ref{sec:superexpanders} for the terminology. By \cite[Theorem 1.3]{MR2316269}, for a locally compact, second countable group, property F$_E$ implies property (T$_E$). The converse is not true: For example, the rank one simple Lie group $\Sp(n,1)$, with $n \geq 2$, has (T$_{L^p}$) for every $p>2$ (as every locally compact, second countable group with property (T) \cite[Theorem A]{MR2316269}), but it has proper actions by isometries on $L^p$-spaces for $p$ large enough \cite{MR1086210,MR2421319} (and so does every hyperbolic group \cite{MR2221161}). In fact, no property (T) group is known to fail property (T$_E$) for some super-reflexive Banach space $E$. On the other hand, the only groups (besides compact groups) known to have property (F$_E$) for every super-reflexive Banach space $E$ are the ones from Theorem~\ref{thm:BFGMconjecture}.

Oppenheim's main result is that $\SL_3(\Z)$ has property (T$_E$) for every super-reflexive Banach space $E$. By previously known arguments, this implies that $\SL_3(\R)$ has (T$_E$), and then that $\SL_{n\geq 4}(\R)$ has F$_E$. More precisely, Oppenheim proves that the pair $(\SL_3(\Z),X_{ij}(\Z))$ has a form of relative property (T$_E$) for every $1 \leq i \neq j\leq 3$, where $X_{ij}(\Z)$ is the elementary subgroup $\{ 1+k e_{i,j} \mid k \in \Z\}$; by a bounded generation argument the main result follows. He does so by constructing a sequence $(\mu_n)$ of almost $X_{ij}(\Z)$-invariant probability measures on $\SL_3(\Z)$ such that $\pi(\mu_n)$ converges in the norm topology of $B(E)$ for every super-reflexive Banach space $E$.

One first difference between our proof and Oppenheim's is that we work directly on the level of the Lie groups $\SL_3(\R)$ and $\Sp_4(\R)$. We could probably make the argument work for the lattices $\SL_3(\Z)$ and $\Sp_4(\Z)$, but this would make the arguments technically more involved. The main challenge and difference is to work directly with actions by affine isometries.

Let us first explain the idea for $G=\SL_3(\R)$. The main result is that there is a net $(\mu_\lambda)$ of probability measures such that whenever $G\curvearrowright E$ is an isometric action on a uniformly convex Banach space $E$, averaging along any fixed orbit with respect to $\mu_\lambda$ gives rise to a net $(\xi_\lambda)$ of vectors in $E$ that converges to a vector that is fixed by $G$. By a version of the Mautner phenomenon, it is even enough to show that this limit vector is fixed by some non-compact subgroup of $G$. The measures $\mu_\lambda$ are not the uniform measures on $K$-double cosets (and are not even $K$-bi-invariant) as usual (see e.g.~\cite{MR2423763,MR1760627,MR2186253,MR1151617}). Instead, the measures are constructed as convolution products of Gaussian measures on elementary subgroups, in a precise order given by the root system $A_2$ (see Figure~\ref{picture:A2}). These measures are related to Oppenheim's measures, but with Gaussian measures replacing uniform measures on integer segments. To show that the net $(\xi_\lambda)$ converges, we prove that it satisfies the Cauchy criterion, so we need to compare the vectors $\xi_\lambda$ for different values of $\lambda$. The main original ingredient we introduce is a non-trivial result about isometric actions on uniformly convex spaces of the Heisenberg group $H_3(\R)$, stated in Proposition~\ref{prop:affine_actions_of_H3R}. More precisely, it bounds the distance between the averages along $H_3(\R)$-orbits with respect to two different probability measures constructed as convolution product of Gaussian measures on $1$-parameter subgroups of $H_3(\R)$, in terms of the growth of the orbit. An interesting aspect of Proposition~\ref{prop:affine_actions_of_H3R} is that the statement becomes essentially void for actions that grow polynomially on the center of $H_3(\R)$. The result becomes much stronger for actions that grow at most logarithmically. Fortunately, thanks to the well-known exponential distorion of unipotents (see Lemma~\ref{lem:orbitgrowth}), this is the case for restrictions of actions of $\SL_3(\R)$. Exploiting the various embeddings of the Heisenberg group in $H_3(\R)$, we can then relate the vectors $\xi_\lambda$ for many different values of $\lambda$. Another crucial but much easier tool we need is a quantitative form of amenability of $H_3(\R)$; see Lemma~\ref{lem:TV_of_gaussians_H3_affine}.

For $G=\Sp_4(\R)$, the general strategy is the same. Here again, the main analytic part happens in the nilpotent part of the Iwasawa decompostion of $G$. This nilpotent group $H$ studied in subsection~\ref{subsection:H} is not as simple as the Heisenberg group. For example, it is three-step nilpotent, whereas $H_3(\R)$ is two-step nilpotent. The facts that the commutator subgroup of $H$ is not central and that the root diagram $C_2$ underlying $\Sp_4(\R)$ has fewer symmetries than the root system $A_2$ (there are two types of roots, namely long and short ones) create some obstacles in the analysis. We are nonetheless able to show similar results for affine actions of $H$ on uniformly convex Banach spaces in Propositions~\ref{prop:changec_H} and \ref{prop:nunutildechangeb}. Subsequently, these local estimates are exploited to show that the convolution products of Gaussian measures on the root groups ordered as in the root diagram $C_2$ (see Figure~\ref{picture:C2}) are Cauchy, and therefore converge, to a point that is necessarily a fixed point. The combinatorics of this step are also more involved than for $\SL_3(\R)$.

\subsection{Organization} Section~\ref{sec:preliminaries} contains some preliminary discussions on various known facts, and also some computations on optimal transport distances between Gaussian measures. Section~\ref{sec:R} contains our first important results about isometric representations of the group $\R$ on uniformly convex Banach spaces. The proof of Theorem~\ref{thm:higherrankgroupsFE} for the cases $\SL_3(\R)$ and $\Sp_4(\R)$ is given in Section~\ref{sec:SL3} and \ref{sec:Sp4}, respectively. Theorem~\ref{thm:higherrankgroupsFE} and Theorem~\ref{thm:BFGMconjecture} are then deduced in Section~\ref{sec:conclusion}. Finally, we explain how to deduce Corollary~\ref{cor:superexpanders} in Section~\ref{sec:superexpanders}.

\section{Preliminaries}\label{sec:preliminaries}

\subsection{Affine isometric actions}
Throughout this article, every action (and in particular every representation) of a topological group $G$ on a topological space $E$ is assumed to be continuous, i.e.~the map $G\times E\to E$ is continuous. Every Banach space is assumed to be defined over the field $\R$; all results will automatically hold for Banach spaces over $\C$ as well.

We review a few standard facts on affine isometric actions and refer to \cite[Section 2.d]{MR2316269} for details. By the Banach--Mazur theorem, isometries of real Banach spaces $E$ are affine. Therefore, if $G\curvearrowright E$ is an action by isometries, it has the form
\[
    g\cdot \xi = \pi(g) \xi+ b(g),
\]
where $\pi\colon G\to \mathrm{O}(E)$ is a continuous isometric representation and $b\colon G\to E$ is a continuous map satisfying the $1$-cocycle relation $b(gh) = b(g) + \pi(g) b(h)$ for all $g,h \in G$.

It is natural to extend any affine isometric action to the space of probability measures as follows. Let $G$ be a compactly generated locally compact group with compact generating set $S$ and associated word-length $|\cdot|_S$. If $G$ acts continuously by affine isometries on a Banach space $E$, and if $\mu$ is any Borel probability measure on $G$ such that $\int |g|_S \, d\mu(g) <\infty$, then for $\xi \in E$, we write $\mu \cdot \xi = \int g \cdot \xi \, d\mu(g)$ (Bochner integral of the continuous map $g\mapsto g\cdot \xi$). The convergence of the integral is justified by the cocycle inequality
\begin{equation} \label{eq:cocycle_inequality}
    \|g \cdot \xi-\xi\| \leq |g|_S \max_{s \in S} \| s \cdot \xi-\xi\|,
\end{equation}
which is a consequence of the subadditivity of the function $g \mapsto \|g\cdot \xi - \xi\|$. It is worth pointing out that this inequality not only justifies the definition of $\mu \cdot \xi$, but also provides us with the estimate
\begin{equation} \label{eq:displacement_affine}
    \|\mu \cdot \xi - \xi\| \leq \int |g|_S  \max_{s \in S}\|s\cdot \xi - \xi\| \, d\mu(g).
\end{equation}

More generally, if $\mu$ and $\nu$ are two such probability measures and $\Pi$ is a coupling between them, i.e.~$\Pi$ is a Borel probability measure on $G \times G$ with first marginal $\mu$ and second marginal $\nu$, then
    \begin{align*}
    \|\mu \cdot \xi - \nu \cdot \xi\| & \leq \int \|g\cdot \xi - h\cdot \xi\| \, d\Pi(g,h)\\
    &= \int \|h^{-1} g\cdot \xi - \xi\| \, d\Pi(g,h) \\
    &\leq \int |h^{-1}g|_S \, d\Pi(g,h) \, \max_{s \in S}\|s\cdot \xi - \xi\|,
    \end{align*}
where the last inequality follows from \eqref{eq:cocycle_inequality}. Thus, if we let $\mathcal{T}_{d_S}(\mu,\nu)$ denote the Kantorovich--Rubinstein (or Wasserstein) distance with cost function $d_S(g,h) = |h^{-1}g|_S$ (see e.g.~\cite[Chapter 6]{MR2459454}), i.e., $\mathcal{T}_{d_S}(\mu,\nu)$ is the infimum of $\int |h^{-1}g|_S \, d\Pi(g,h)$ over all couplings $\Pi$ of $\mu$ and $\nu$, we deduce 
\begin{equation} \label{eq:cocycle_and_transport}
    \|\mu \cdot \xi - \nu \cdot \xi\| \leq \mathcal{T}_{d_S}(\mu,\nu) \, \max_{ s \in S}\|s\cdot \xi - \xi\|.
\end{equation}

\subsection{Gaussian measures}
For $t \in \R$ and $a \in \R$, let $\gamma_{t,a}$ denote the Gaussian probability measure with mean $t$ and variance $e^{2a}$. For the centered Gaussian probability measure (i.e.~the case $t=0$) with variance $e^{2a}$, we write $\gamma_a$. We will use this notation throughout this article.

We recall the following standard fact, which is obtained by straightforward computation. In this lemma and in the rest of this article, we use the normalization of the total variation distance that gives values between $0$ and $1$, that is, the total variation distance between two probability measures on $(X,\mathcal B)$ is given by
\[
    \|\mu - \nu\|_{TV}  = \sup\{ |\mu(A) - \nu(A)|  \mid A \in \mathcal B\} \in [0,1].
\]

\begin{lemma} \label{lem:TV_of_gaussians}
For all $t,a \in \R$, the total variation distance between $\gamma_{t,a}$ and $\gamma_a$ is bounded above by $\frac{|t|}{\sqrt{2\pi}}e^{-a}$.
\end{lemma}

The total variation distance coincides with the optimal transportation cost for the cost function $1_{x\neq y}$. We also need the similar fact for the cost function $c(x,y) = |x-y|+1_{x \neq y}$. Following the previous notation, if $\mu,\nu$ are probability measures on $\R$, let 
    \[
        \mathcal{T}_c(\mu,\nu) = \inf \int (|x-y|+1_{x \neq y}) \, d\Pi(x,y)
    \]
    denote the Kantorovich--Rubinstein distance associated with $c$, where the infimum is taken over all probability measures $\Pi$ on $\R^2$ with first marginal $\mu$ and second marginal $\nu$.
\begin{lemma} \label{lem:Monge_problem}
For all $s,t \in \R$ and $b\geq 0$, we have $\mathcal{T}_c(\gamma_{s,0},\gamma_{t,b}) \leq 2(|s-t|+e^b-1)$.
%\footnote{The proof gives the following more precise bound, valid for all real numbers $a,b,s,t$ with $a\leq b$: \[ \mathcal{T}_c(\gamma_{s,a},\gamma_{t,b}) \leq (\frac{1}{\sqrt{2\pi}} e^{-a} +1)|t-s| + \frac{2 e^{(-1/2)}}{\sqrt{2\pi}}(b-a) + \frac{1}{\sqrt{2\pi}}(8 e^{(-1/2)}-2)(e^b - e^a).\]}
\end{lemma}
\begin{proof}
Since $\mathcal{T}_c$ satisfies the triangle inequality (see \cite[Chapter 6]{MR2459454}), it suffices to prove the following two inequalities for every $s \in \R$ and $a \geq 0$:
\begin{equation} \label{eq:Monge_mean}
    \limsup_{t\to s} \frac{1}{|t-s|} \, \mathcal{T}_c(\gamma_{s,0},\gamma_{t,0}) \leq 2,
\end{equation}
\begin{equation} \label{eq:Monge_variance}
    \limsup_{b\to a}\frac{1}{|e^b-e^a|} \, \mathcal{T}_c(\gamma_{s,a},\gamma_{s,b}) \leq 2.
\end{equation}
Since the cost $|x-y|+1_{x \neq y}$ only depends on $|x-y|$, we can assume $s=0$.

Before we show \eqref{eq:Monge_mean} and \eqref{eq:Monge_variance}, we show an obvious (but useful) general upper bound for $\mathcal{T}_c(\mu,\nu)$ for two absolutely continuous probability measures $d\mu = f \, dx$ and $d\nu = g \, dx$. The measures $(f-g)_+ \, dx$ and $(g-f)_+ \, dx$ have the same total mass. Let $\Pi_0$ be any coupling between them, and define a measure $\Pi$ on $\R^2$ by
\[
    \int h(x,y) \, d\Pi(x,y) = \int h(x,x) \min(f(x), g(x)) \, dx + \int h(x,y) \, d\Pi_0(x,y).
\]
The measure $\Pi$ is a coupling between $\mu$ and $\nu$, so
\begin{align*}
    \mathcal{T}_c(\mu,\nu) &\leq \int (|x-y|+1_{x \neq y}) \, d\Pi(x,y) \\
    &= \int (1+|x-y|) \, d\Pi_0(x,y)\\ 
    & \leq \int \left(\frac 1 2 +|x|+\frac 1 2+|y|\right) \, d\Pi_0(x,y)\\ 
    & = \int \left(\frac{1}{2}+|x|\right) |f(x) - g(x)| \, dx.
\end{align*}  
The third line is the obvious inequality $|x-y| \leq |x|+|y|$, and the fourth line is because the marginals of $\Pi_0$ are $(f-g)_+ dx$ and $(g-f)_+ dx$. Going back to \eqref{eq:Monge_mean} and \eqref{eq:Monge_variance}, we see that it suffices to prove that,
\begin{align*}
    & \limsup_{t\to 0} \int \left(\frac{1}{2}+|x|\right) \frac{|f_{t,0}(x) - f_{0,0}(x)|}{|t|} \, dx \leq 2, \\
    \forall a\geq 0, \qquad & \limsup_{b\to a} \int \left(\frac{1}{2}+|x|\right) \frac{|f_{0,b}(x) - f_{0,a}(x)|}{|e^b - e^a|} \, dx \leq 2,
\end{align*}
where $f_{s,a}(x) =  \frac{1}{\sqrt{2\pi}e^a} e^{-\frac{(x-s)^2}{2e^{2a}}}$ is the density of $\gamma_{s,a}$. We can safely exchange the limit and the integral, so we are left to prove that
\begin{align*}
    & \int \left(\frac{1}{2}+|x|\right) |x| \, f_{0,0}(x) \, dx \leq 2, \\
    \forall a\geq 0, \qquad &\int \left(\frac{1}{2}+|x|\right) e^{-a} \, |1-x^2e^{-2a}| \, f_{0,a}(x) \, dx \leq 2.
\end{align*}
The first integral is equal to $1+\frac{1}{\sqrt{2\pi}}$, which is indeed less than $2$. By the change of variable $u=xe^{-a}$, the second inequality becomes
\[
    \int \left(\frac{1}{2e^a}+|u|\right) |1-u^2| \, f_{0,0}(u) \, du \leq 2.
\]
Since this expression is decreasing in $a$ and $a\geq 0$, we are left to observe that
\[
    \frac{1}{\sqrt{2\pi}} \int_\R \left(\frac 1 2+|u|\right) |1-u^2| e^{-\frac{u^2}{2}} \, du\leq 2.
\]
By explicit computation, the above integral is equal to $\frac{1}{\sqrt{2\pi}}\left(\frac{10}{\sqrt{e}}-2\right) \simeq 1.6218$.
\end{proof}

\subsection{Orbit growth of actions by unipotent elements}
We will often use the following fact, which combines exponential distortion of unipotent elements in a semisimple Lie group (see e.g.~\cite[Section 3]{MR1253544}) with linear orbit growth of affine isometric actions. We state it for semisimple Lie groups (and only use it for the groups $\SL_3(\R)$ and $\Sp_4(\R)$), but it holds much more generally.

\begin{lemma} \label{lem:orbitgrowth}
    Let $G$ be a semisimple Lie group with Lie algebra $\mathfrak{g}$, and let $X \in \mathfrak{g}$ be a nilpotent element. For every action of $G$ by affine isometries on a Banach space $E$ and every $\xi \in E$, there is a constant $M$ (depending on $X$ and $\xi$) such that for every $t\in \R$,
    \[
        \|\exp(tX) \cdot \xi - \xi\| \leq M \log(2+|t|).
    \]
\end{lemma}
\begin{proof}
We can assume $t \geq 1$, because the inequality is obvious for $|t|\leq 1$ and we have $\| \exp(-tX) \cdot \xi - \xi\| = \|\xi -  \exp(tX) \cdot \xi \|$.

By the Jacobson--Morozov theorem, there exists $Y \in \mathfrak g$ such that $[Y,X]=X$. Taking the exponential, we deduce
    \[
        \exp(sY) \exp(X) \exp(-sY) = \exp( e^s X).
    \]
    Taking $s=\log t$ and applying the cocycle inequality~\eqref{eq:cocycle_inequality}, we deduce
    \begin{align*}
        \| \exp(tX) \cdot \xi - \xi\| & = \| \exp(s Y) \exp(X) \exp(-sY) \cdot \xi - \xi\|\\
        &\leq\|\exp(X)\cdot \xi - \xi\|+2(1+|s|)\max_{0\leq u \leq 1} \|\exp(uY) \cdot \xi - \xi\|. 
    \end{align*}
This proves the lemma.
\end{proof}

\section{Actions of $\R$}\label{sec:R}

This section is devoted to affine actions of $\R$ on Banach spaces. The main results in this section, Proposition~\ref{prop:unif_convexityR} and \ref{prop:uniform_convexity_iteratedR}, deal with uniformly convex spaces. We start with some easier facts, which hold without restriction on the Banach space.

\begin{lemma} \label{lem:TV_of_gaussiansbis}
    Let $\R$ act continuously by affine isometries on a Banach space $E$, and let $\xi \in E$. For all real numbers $s,t,a,b$ with $a\leq b$,
\begin{equation} \label{eq:tail_of_gaussians}
    \| \gamma_a \cdot \xi - \xi\| \leq \left(1+\sqrt{\frac{2}{\pi}}\right) \, \max_{|u|\leq e^{a}}\|u\cdot \xi - \xi \|.
\end{equation} 
and
\begin{equation} \label{eq:actions_of_gaussians}
            \| \gamma_{s,a}\cdot \xi - \gamma_{t,b} \cdot \xi\| \leq \left(2 \frac{|s-t|}{e^a}+\frac{|e^{2b} - e^{2a}|}{e^{2a}}\right) \, \max_{|u|\leq e^a} \|u\cdot \xi-\xi\|.
\end{equation}
\end{lemma}
\begin{proof}
    By homogeneity we can assume $a=0$. In that case the first inequality was already proved in \eqref{eq:displacement_affine}, because the word-length with respect to the generating set $[-1,1]$ is bounded above by $1+|t|$ and $\int |t| \, d\gamma_0(t) = \sqrt{\frac 2 \pi}$. The second inequality is immediate from \eqref{eq:cocycle_and_transport} and Lemma~\ref{lem:Monge_problem}, using that the word-length on $\R$ is bounded above by the cost function $c$, and using $2(e^x-1) \leq e^{2x}-1$.
\end{proof}

\subsection{Uniformly convex Banach spaces} \label{subsec:ucsr}
A Banach space $E$ is uniformly convex if for every $\varepsilon \in (0,2]$, there exists $\delta > 0$ such that for all $\xi,\eta \in E$ with $\|\xi\|=\|\eta\|=1$,
\[
    \|\xi - \eta\| \geq \varepsilon \quad \Longrightarrow \quad \left\| \frac{\xi+\eta}{2} \right\| \leq 1-\delta.
\]
A Banach space $E$ is super-reflexive if every Banach space that is finitely representable in $E$ is reflexive. Every uniformly convex Banach space is super-reflexive, and every super-reflexive Banach space admits an equivalent uniformly convex norm \cite{MR0336297}.

We state the following elementary lemma for future reference.
\begin{lemma} \label{lem:from_pointwise_to_gaussian}
        Let $E$ be a Banach space and $\rho \colon \R \to \O(E)$ an isometric \hyphenation{rep-re-sen-ta-tion}representation. For all $a \in \R$ and $\xi \in E$,
       \[
            \| \rho(\gamma_a) \xi-\xi\| \leq \left(1+\sqrt{\frac 2 \pi} \right) \max_{|t|\leq e^a} \|\rho(t) \xi - \xi\|.
        \]
\end{lemma}
\begin{proof}
This is \eqref{eq:tail_of_gaussians} for the action $t \cdot \xi = \rho(t) \xi$. 
\end{proof}
In particular, if $\max_{|t|\leq e^a} \|\rho(t) \xi - \xi\|$ is much smaller than $\|\xi\|$, then $\|\rho(\gamma_a)\xi\|$ is close to $\|\xi\|$. The following proposition provides a form of the converse, under the assumption that $E$ is uniformly convex. This is a form of the general and well-understood phenomenon that for representations on uniformly convex spaces, having almost invariant vectors is equivalent to the averaging operator with respect to reasonable probability measures having norm $1$ (see \cite[Proposition 5.1]{MR3781331} or \cite[Theorem 3.4]{MR4000570}). We will need a quantitative form, which is obtained by the same proof as in \cite[Theorem 3.4]{MR4000570}.
\begin{prop}\label{prop:unif_convexityR} Let $E$ be a uniformly convex Banach space. For every $\varepsilon>0$, there is $\delta>0$ such that the following holds: For all isometric representations $\rho \colon \R \to \O(E)$, all $a \in \R$, and all $\xi \in E$,
\[
    \| \rho(\gamma_a) \xi \|\geq (1-\delta) \|\xi\| \implies \max_{ |t| \leq e^a }\|\rho(t) \xi - \xi\| \leq \varepsilon\|\xi\|.
\]
\end{prop}
\begin{proof} By homogeneity, we can assume that $a=0$. Let $\varepsilon > 0$. By the uniform convexity of $E$, there is $\delta_0$ such that for every two unit vectors $\xi,\eta \in E$,
\[
    \left\|\frac{\xi+\eta}{2}\right\| \geq 1-\delta_0 \implies \|\xi-\eta\| \leq \frac \varepsilon 2.
\]
Let $c<1$ be the total variation distance between $\gamma_0$ and $\gamma_{1,0}$. We prove the proposition with $\delta = (1-c) \delta_0$.  
  
For $t \in (-1,1)$, the total variation distance between $\gamma_0$ and $\gamma_{t,0}$ is less than $c$. Let $(X,X')$ be a coupling achieving the total variation distance, i.e.~$X$ is distributed as $\gamma_0$, $X'$ is distributed as $\gamma_{t,0}$, and $\|\gamma_0-\gamma_{t,0}\|_{TV} = \prob(X \neq X') \leq c$. Then
\[
    \rho(\gamma_0) \xi = \E \rho(X) \xi = \E \rho(X'-t)\xi = \E \frac{\rho(X) + \rho(X'-t)}{2}\xi.
\]
If $X=X'$, then $\rho(X)\xi + \rho(X'-t)\xi = \rho(X-t)( \rho(t) \xi+\xi)$ has the same norm as $\rho(t) \xi + \xi$. Otherwise, we can estimate the norm by $2\|\xi\|$ by the triangle inequality. Therefore, if $\| \rho(\gamma_0) \xi \|\geq (1-\delta) \|\xi\|$, we have
\begin{align*}
    (1-\delta)\|\xi\| & \leq \E\frac 1 2 \|(\rho(X) + \rho(X'-t))\xi\|\\
& \leq \prob(X \neq X')\|\xi\| + \prob(X = X') \frac 1 2 \|\rho(t) \xi + \xi\|\\
& \leq c\|\xi\|+(1-c) \frac 1 2 \|\rho(t) \xi + \xi\|.
\end{align*}
Equivalently, 
\[
    \left(1-\frac{\delta}{1-c} \right) \|\xi\| \leq \frac 1 2 \|\rho(t) \xi + \xi\|.
\]
By our choice of $\delta=(1-c)\delta_0$ and of $\delta_0$, we obtain $\|\rho(t) \xi - \xi\| \leq \varepsilon \|\xi\|$.
\end{proof}

The following proposition generalises the result above to the case when we have several representations of $\R$. This will be used in a situation where these representations are restrictions of a representation of single group to various one-parameter subgroups. It is worth noting that Proposition~\ref{prop:uniform_convexity_iteratedR} cannot be obtained as a consequence of \cite[Theorem 3.4]{MR4000570}, because the convolution of Gaussian measures on different one-parameter subgroups is not at all reasonable. Forms of this result appear somewhat implicitely in the proof of \cite[Proposition 5.4]{Oppenheim2022}.

\begin{prop} \label{prop:uniform_convexity_iteratedR} 
Let $E$ be a uniformly convex Banach space and $n$ a positive integer. For every $\varepsilon>0$, there is $\delta>0$ such that the following holds: For all isometric representations $\rho_1,\dots,\rho_n \colon \R \to \O(E)$, all $a_1,\dots a_n \in \R$, and all $\xi \in E$,
\begin{equation} \label{eq:uniform_convexity_iterated}
    \|\rho_1(\gamma_{a_1}) \dots \rho_n(\gamma_{a_n}) \xi \|\geq (1-\delta) \|\xi\| \implies\max_{1\leq i \leq n} \max_{|t| \leq e^{a_i}} \|\rho_i(t) \xi - \xi\| \leq \varepsilon \|\xi\|.
\end{equation}
\end{prop}
\begin{proof}
By homogeneity we can assume that $a_i=0$ for all $i$. 

By induction on $n$, we can prove the stronger implication
\begin{equation} \label{eq:uniform_convexity_iterated_2}
    \|\rho_1(\gamma_{0}) \dots \rho_n(\gamma_{0}) \xi \|\geq (1-\delta)\|\xi\| \implies \begin{cases} \max_i \max_{|t| \leq 1} \|\rho_i(t) \xi - \xi\| \leq \varepsilon\|\xi\|, \\ \max_i \|\rho_i(\gamma_{0}) \xi - \xi\| \leq \varepsilon\|\xi\|.\end{cases}
\end{equation}
Indeed, let $\varepsilon_n \colon [0,1]\to[0,2]$ be the function defined by setting $\varepsilon_n(\delta)$ to be the smallest $\varepsilon$ such that \eqref{eq:uniform_convexity_iterated_2} holds for all $n$, $\rho_1,\dots,\rho_n$ and $\xi$. The assertion of the proposition exactly says that \begin{equation} \label{eq:epsilon_nC0}
    \lim_{\delta \to 0}\varepsilon_{n}(\delta)=0.
\end{equation}
The case $n=1$ follows from Proposition~\ref{prop:unif_convexityR} (with a different $\delta$), using Lemma~\ref{lem:from_pointwise_to_gaussian}.

Now, suppose that \eqref{eq:uniform_convexity_iterated_2} holds for all $k<n$, let $\rho_1,\dots,\rho_n$ and $\xi$ be given, and moreover, suppose that
\[
    \|\rho_1(\gamma_{0}) \dots \rho_n(\gamma_{0}) \xi \|\geq (1-\delta)\|\xi\|.
\]
In particular, we have $\|\rho_2(\gamma_{0}) \dots \rho_n(\gamma_{0}) \xi \|\geq (1-\delta)\|\xi\|$, and therefore,
\[
    \begin{cases}
        \max_{2\leq i \leq n} \max_{|t| \leq 1} \|\rho_i(t) \xi - \xi\| \leq \varepsilon_{n-1}(\delta)\|\xi\|, \\ \max_{2\leq i \leq n} \|\rho_i(\gamma_{0}) \xi - \xi\| \leq \varepsilon_{n-1}(\delta)\|\xi\|.
    \end{cases}
\]
By the triangle inequality, $\| \rho_1(\gamma_0) \dots \rho_n(\gamma_0) \xi - \rho_1(\gamma_0) \xi\| \leq (n-1) \varepsilon_{n-1}(\delta)\|\xi\|$, which implies $\|\rho_1(\gamma_0)\xi\| \geq (1-\delta - (n-1) \varepsilon_{n-1}(\delta)) \|\xi\|$, and hence,
\begin{align*}
    \max_{|t|\leq 1} \|\rho_1(t) \xi - \xi\| &\leq \varepsilon_1(\delta+(n-1) \varepsilon_{n-1}(\delta)) \|\xi\| \quad \textrm{and} \\
    \|\rho_1(\gamma_0) \xi - \xi\| &\leq \varepsilon_1(\delta+(n-1) \varepsilon_{n-1}(\delta)) \|\xi\|. 
\end{align*}
Altogether, we have shown that
\[
    \varepsilon_n(\delta) \leq \max(\varepsilon_{n-1}(\delta), \varepsilon_1(\delta+(n-1) \varepsilon_{n-1}(\delta))).
\]
\end{proof}

\section{The group $\SL_3(\R)$} \label{sec:SL3}

In this section, we prove Theorem \ref{thm:higherrankgroupsFE} for $G=\SL_3(\R)$. As explained in Section \ref{subsec:strategy}, we use the same structure theory as Oppenheim in \cite{Oppenheim2022}, but we directly consider isometric actions rather than isometric representations.

\subsection{Representations and actions of the real Heisenberg group $H_3(\R)$}\label{subsec:heisenberg}
In this section, we fix a uniformly convex Banach space $E$.

Consider the Heisenberg Lie algebra $\mathfrak{h}_3(\R)$, i.e.~the $3$-dimensional Lie algebra with basis $\mathfrak{X},\mathfrak{Y},\mathfrak{Z}$ and relations
\[
	[\mathfrak{X},\mathfrak{Y}]=\mathfrak{Z}, \qquad [\mathfrak{X},\mathfrak{Z}]=[\mathfrak{Y},\mathfrak{Z}]=0.
\]
Let $H_3(\R)$ be the simply connected Lie group with Lie algebra $\mathfrak{h}_3(\R)$, and let $\exp \colon \mathfrak{h}_3(\R) \to H_3(\R)$ denote the associated exponential map. For $r,s,t \in \R$, the elements $X(t)=\exp(t\mathfrak{X})$, $Y(s) = \exp(s\mathfrak{Y})$, and $Z(r)=\exp(r\mathfrak{Z})$ of $H_3(\R)$ satisfy
\[
    [X(t),Y(s)]=X(t)^{-1}Y(s)^{-1}X(t)Y(s)=Z(ts).
\]

We define two families of probablity measures on $H_3(\R)$ constructed as convolution products of Gaussian measures on the groups $X(\R)$, $Y(\R)$, and $Z(\R)$. Specifically, for real numbers $a,b,c$, let $\nu_{abc}$ and $\widetilde{\nu}_{abc}$ be the probability measures on $H_3(\R)$ defined as
\[
    \nu_{abc} = X(\gamma_a) Z(\gamma_{b}) Y(\gamma_{c}) \quad \textrm{and} \quad \widetilde{\nu}_{abc} = Y(\gamma_{a}) Z(\gamma_{b}) X(\gamma_{c}).
\]

An easy application of Lemma~\ref{lem:TV_of_gaussians} gives the following.
\begin{lemma}\label{lem:TV_of_gaussians_H3}
For all $a,b,c \in \R$,
    \[
        \|\nu_{abc} - \widetilde{\nu}_{cba}\|_{TV} \leq e^{a+c-b}.
    \]
\end{lemma}
\begin{proof}
    Using the fact that $Z$ takes values in the center, we can write $\nu_{abc} - \widetilde{\nu}_{cba}$ as the average, when $t\sim \gamma_{a}$ and $s \sim \gamma_{c}$, of the measures $X(t) Y(s) \cdot Z(\gamma_{b}) - Y(s) X(t)  \cdot Z(\gamma_{b})$, where $\cdot$ denotes the action of $H_3(\R)$ by left-translation on the measures. By the triangle inequality, we deduce
\begin{align*}
    \|\nu_{abc} - \widetilde{\nu}_{cba}\|_{TV} &\leq \E_{t \sim \gamma_{a}, s \sim \gamma_{c}} \|X(t)Y(s)\cdot Z(\gamma_{b})-Y(s)X(t)\cdot Z(\gamma_{b})\|_{TV} \\
    &= \E_{t \sim \gamma_{a}, s \sim \gamma_{c}} \|X(t)^{-1}Y(s)^{-1}X(t)Y(s)\cdot Z(\gamma_{b})-Z(\gamma_{b})\|_{TV} \\
    &= \E_{t \sim \gamma_{a}, s \sim \gamma_{c}} \|Z(ts)\cdot Z(\gamma_{b})-Z(\gamma_{b})\|_{TV} \\
    &= \E_{t \sim \gamma_{a}, s \sim \gamma_{c}} \|Z(\gamma_{ts,b})-Z(\gamma_{b})\|_{TV} \\
    & \leq \E_{t \sim \gamma_{a}, s \sim \gamma_{c}} \frac{|t||s|}{\sqrt{2\pi}} e^{-b}\\
    & = \frac{\sqrt 2}{\pi\sqrt{\pi}} e^{a+c-b} \leq e^{a+c-b}.
\end{align*}
The last equality is just the fact that $\E_{t\sim \gamma_{a}} |t| = \sqrt{\frac 2 \pi} e^{a}$.
\end{proof}

\begin{lemma} \label{lem:pi(X0Yb)_contractionR}
  There exists $q_0 < 1$ such that for all isometric representations $\pi \colon H_3(\R) \to \O(E)$, all $a,c \in \R$, and all $\xi \in E$,
\begin{equation} \label{eq:pi(X0Yb)_contraction}
    \|\pi(X(\gamma_{a}) Y(\gamma_{c}))\xi\| \leq \max\left( q_0 \|\xi\|, 2 \|\pi(Z(\gamma_{a+c})) \xi\|\right).
\end{equation}
\end{lemma}
\begin{proof}
Let $\varepsilon = \frac{1}{8\left(1+\sqrt{\frac 2 \pi}\right)}$, and let $\delta>0$ be given by Proposition~\ref{prop:uniform_convexity_iteratedR} for $n=2$ and this $\varepsilon$. We prove the lemma for $q_0=1-\delta$. Let $\pi$, $a$, $c$, and $\xi$ be given. Suppose $\|\xi\|=1$ and $\|\pi(X(\gamma_{a})Y(\gamma_{c}))\xi\| \geq q_0$. By Proposition~\ref{prop:uniform_convexity_iteratedR}, $\|\pi(X(t)) \xi - \xi\| \leq \varepsilon$ for every $|t|\leq e^{a}$ and $\|\pi(Y(s))\xi - \xi\| \leq \varepsilon $ for every $|s| \leq e^{c}$. This implies that 
\[
    \| \pi(Z(ts)) \xi - \xi\| = \|\pi([X(t),Y(s)])\xi - \xi\| \leq 4 \varepsilon.
\]
In other words, $\|\pi(Z(r))\xi-\xi\| \leq 4\varepsilon$ for every $|r| \leq e^{a+c}$ and therefore, by Lemma~\ref{lem:from_pointwise_to_gaussian},
\[
    \|\pi(Z(\gamma_{a+c}))\xi - \xi\| \leq 4\varepsilon \left( 1+\sqrt{\frac 2 \pi} \right) = \frac 1 2.
\]
As a consequence, $\|\pi(Z(\gamma_{a+c}))\xi\| \geq \frac 1 2$, which proves the lemma.
\end{proof}

\begin{prop} \label{prop:affine_actions_of_H3R}
There exist $q < 1$ and $C > 0$ such that for all isometric actions $H_3(\R) \curvearrowright E$, all $a,b,b',c \in \R$ with $\max(b,b') \leq a+c$, and all $\xi \in E$,
\[
    \| \nu_{abc} \cdot \xi - \nu_{ab'c} \cdot \xi\| \leq C q^{\Delta} (1+|b-b'|)\max_{|t| \leq e^{a+c}} \|Z(t) \cdot \xi-\xi\|,
\]
where
\[
    \Delta=\sqrt{a+c-\max(b,b')}.
\]
\end{prop}

In order to prove the proposition, we use the following lemma.

\begin{lemma} \label{lem:calculus}
For all $a,b,c$ and $\Delta$ as in Proposition~\ref{prop:affine_actions_of_H3R} with $\Delta \geq 3$, there exists an integer $n$ with $\frac{\Delta}{9} \leq n \leq \Delta$ and real numbers $a_1,\dots,a_n$, $c_1,\dots,c_n$ such that
\begin{enumerate}[(i)]
    \item\label{item1} $a_i+c_i \geq b+\Delta$ for all $i$,
    \item\label{item2} $a_i + c_j \leq b-\Delta$ for all $i>j$,
    \item\label{item3} $\sum_{i=1}^n e^{2a_i} = e^{2a}$ and $\sum_{i=1}^n e^{2c_i} = e^{2c}$.
\end{enumerate}
\end{lemma}
\begin{proof}
Set $n=\lfloor \frac{\Delta}{9} \rfloor+1$, $a_i = a-(2i-1) \Delta$ and $c_i = b+\Delta-a_i$ for $i = 1, \ldots, n$. It is a small computation to verify (\ref{item1}) and (\ref{item2}) and to show $\sum_{i=1}^n e^{2a_i} \leq e^{2a}$. Furthermore, observe that $c_i \leq c_n$ for $i=1,\ldots,n$ and $c_n = b-a+2n\Delta \leq b-a+\frac{2}{3}\Delta^2 \leq c - \frac{\Delta^2}{3}$. It follows that $\sum_{i=1}^n e^{2c_i} \leq \sum_{i=1}^n n e^{2c_n} \leq e^{2c}$. By increasing $a_1$ and $c_n$, which does not affect assertions (\ref{item1}) and (\ref{item2}), we can make assertion (\ref{item3}) hold.
\end{proof}

\begin{proof}[Proof of Proposition~\ref{prop:affine_actions_of_H3R}] If $0 \leq \Delta \leq 3$, by \eqref{eq:tail_of_gaussians}, the proposition is obvious as soon as $C q^3 \geq 2\left(1+\sqrt{\frac 2 \pi}\right)$. We can therefore assume $\Delta \geq 3$. 

Consider an arbitrary isometric action $H_3(\R) \curvearrowright E$. We first make the additional assumption that $b < b' \leq b+1$. Define $b''\leq b'$ by $e^{2b'} = e^{2b} + e^{2b''}$.

Let $n$, $a_1,\ldots,a_n$, and $c_1,\ldots,c_n$ be given by Lemma~\ref{lem:calculus}. For $1 \leq k \leq n$, define $\tilde{a}_k,\tilde{c}_k$ by $e^{2\tilde{a}_k} = \sum_{j=k}^n e^{2a_j}$ and $e^{2\tilde{c}_k} = \sum_{j=k}^n e^{2c_j}$, set
\[
    \eta_k = \nu_{\tilde{a}_k b \tilde{c}_k}\cdot \xi - \nu_{\tilde{a}_k b' \tilde{c}_k}\cdot \xi
\]
and extend this by
\[
    \eta_{n+1} = Z(\gamma_b)\cdot \xi - Z(\gamma_{b'}) \cdot \xi.
\]
By Lemma~\ref{lem:TV_of_gaussiansbis} and the inequalites $b'' \leq b' \leq a+c$, we have
\begin{equation} \label{eq:norm:of_etan+1}
    \|\eta_{n+1}\| \leq  \|Z(\gamma_{b''})\cdot \xi - \xi\| \leq 2 \max_{|t| \leq e^{a+c}} \|Z(t) \cdot \xi - \xi\|.
\end{equation}
For every $k\leq n-1$, using that $\gamma_{\tilde{a}_k} =  \gamma_{a_k} \gamma_{\tilde{a}_{k+1}}$ and $\gamma_{\tilde{c}_k} = \gamma_{\tilde{c}_{k+1}} \gamma_{c_k}$, and denoting the isometric representation underlying the action by $\pi$, we write 
\begin{equation} \label{eq:diff_etak_pi_etak+1}
    \eta_k - \pi(X(\gamma_{a_k}) Y(\gamma_{c_k})) \eta_{k+1} = \pi(X(\gamma_{a_k}) m Y(\gamma_{\tilde{c}_{k+1}}))(\xi - Z(\gamma_{b''})\cdot\xi)
\end{equation}
for the signed measure
\[
    m = \bigg( X(\gamma_{\tilde{a}_{k+1}})Y(\gamma_{c_k}) - Y(\gamma_{c_k})X(\gamma_{\tilde{a}_{k+1}}) \bigg) Z(\gamma_b)  = \nu_{\tilde{a}_{k+1} b c_k} - \widetilde{\nu}_{c_k b \tilde{a}_{k+1}}.
\]
By Lemma~\ref{lem:TV_of_gaussians_H3}, we have
\[
    \|m\|_{TV} \leq e^{\tilde{a}_{k+1} + c_k -b} \leq \sqrt{n} \max_{j\geq k+1} e^{a_j+c_k-b} \leq \sqrt{n} e^{-\Delta}.
\]
The second inequality uses $a_j+c_k-b \leq 0$ and the last inequality uses Lemma~\ref{lem:calculus}.(\ref{item2}).

On the other hand, by \eqref{eq:norm:of_etan+1}, we have
\[
    \|\pi(Y(\gamma_{\tilde{c}_{k+1}})) (\xi - Z(\gamma_{b''})\cdot\xi)\| \leq 2 \max_{|t| \leq e^{a+c}} \|Z(t) \cdot \xi-\xi\|.
\]
From \eqref{eq:diff_etak_pi_etak+1}, we therefore deduce
\begin{equation} \label{eq:diff_etak_pi_etak+1_norm}
    \|\eta_k - \pi(X(\gamma_{a_k}) Y(\gamma_{c_k}))\eta_{k+1}\| \leq 2 \sqrt{n}e^{-\Delta} \max_{|t| \leq e^{a+c}} \|Z(t) \cdot \xi-\xi\|.
\end{equation}
For $k=n$, the left-hand side vanishes, so \eqref{eq:diff_etak_pi_etak+1_norm} remains true.

By Lemma~\ref{lem:pi(X0Yb)_contractionR} and Lemma~\ref{lem:calculus}.(\ref{item1}),
\[
    \|\pi(X(\gamma_{a_k}) Y(\gamma_{c_k})) \eta_{k+1}\| \leq q_0\|\eta_{k+1}\| + 2 \| \pi(Z(\gamma_{b+\Delta})) \eta_{k+1}\|.
\]
Since $Z(\R)$ is contained in the center of $H_3(\R)$,
\[
    \pi(Z(\gamma_{b+\Delta})) \eta_{k+1} = \pi(\nu_{\tilde{a}_{k+1} b \tilde{c}_{k+1}})(Z(\gamma_{b+\Delta})\cdot \xi - Z(\gamma_{b+\Delta})Z(\gamma_{b''})\cdot \xi).
\]
Using \eqref{eq:actions_of_gaussians}, we have
\[
    \|\pi(Z(\gamma_{b+\Delta})) \eta_{k+1} \|\leq e^{-2\Delta} \big( e^{2b'-2b} - 1 \big) \max_{|t| \leq e^{b+\Delta}} \| Z(t) \cdot \xi - \xi\|.
\]
Also, by our assumption $\Delta \geq 3$, we have $b+\Delta \leq b+\Delta^2 \leq a+c$ and also $e^{-2\Delta}(e^{2b'-2b}-1) \leq e^{-\Delta}$. Therefore, from \eqref{eq:diff_etak_pi_etak+1_norm}, we obtain
\[
    \|\eta_k\| \leq q_0 \|\eta_{k+1}\| + 2(\sqrt{n} +1) e^{-\Delta} \max_{|t| \leq e^{a+c}} \|Z(t) \cdot \xi-\xi\|,
\]
and
\[
    \|\nu_{abc}\cdot \xi - \nu_{ab'c}\cdot\xi\| = \|\eta_1\| \lesssim q_0^n\|\eta_{n+1}\| + 2n(\sqrt{n} +1) e^{-\Delta} \max_{|t| \leq e^{a+c}} \|Z(t) \cdot \xi-\xi\|.
\]
By \eqref{eq:norm:of_etan+1} and the inequalities $\frac{\Delta}{9}\leq n \leq \Delta$, we obtain
\[
    \|\nu_{abc}\cdot \xi - \nu_{ab'c}\cdot\xi\| \lesssim q^\Delta \max_{|t| \leq e^{a+c}} \|Z(t) \cdot \xi-\xi\|,
\]
for $q < \max( q_0^{\frac 1 9},e^{-1})$, which does not depend on the chosen action.

Hence, in the case $b < b' \leq b+1$, we obtain the proposition without the factor $1+|b-b'|$. The factor $1+|b-b'|$ appears when dealing with $b,b'$ possibly far apart, by considering a path $b_1,...,b_n$ from $b$ to $b'$ with $|b_i-b_i'| \leq 1$ for $i=1,\ldots,n-1$, applying the inequality we have just obtained to each term $\| \nu_{ab_ic} \cdot \xi - \nu_{ab_{i+1}c} \cdot \xi\|$. The number of such terms is at most $1+|b-b'|$.  
\end{proof}

We  also need the following variant of Lemma~\ref{lem:TV_of_gaussians_H3}.
\begin{lemma} \label{lem:TV_of_gaussians_H3_affine}
For all isometric actions $H_3(\R) \curvearrowright E$, all $a,b,c \in \R$, and all $\xi \in E$,
\[
    \| \nu_{abc} \cdot \xi - \widetilde{\nu}_{cba} \cdot \xi\| \leq 2 e^{a+c-b} \max_{|r| \leq e^{b}}\|Z(r)\cdot \xi-\xi \|.
\]
\end{lemma}
\begin{proof}
    Similar to the proof of Lemma~\ref{lem:TV_of_gaussians_H3}, we obtain
\[
    \|\nu_{abc} \cdot \xi - \widetilde{\nu}_{cba} \cdot \xi\| \leq \E_{t \sim \gamma_{a}, s \sim \gamma_{c}} \|Z(\gamma_{ts,b}) \cdot \xi -Z(\gamma_{b}) \cdot \xi\|,
\]
which is, by Lemma~\ref{lem:TV_of_gaussiansbis}, bounded above by
\begin{align*}
    2 \; \E_{t \sim \gamma_{a}, s \sim \gamma_{c}} \frac{|t||s|}{e^{b}} \max_{|r| \leq e^{b}}\|Z(r)\cdot \xi-\xi \| = \frac{4}{\pi} e^{a+c-b} \max_{|r| \leq e^{b}} \|Z(r)\cdot \xi-\xi \|.
\end{align*}
%If, instead of the crude bound of Lemma~\ref{lem:TV_of_gaussiansbis}, we use the better bound in the footnote of Lemma~\ref{lem:Monge_problem}, we obtain $\frac 2 \pi (1+\frac{1}{\sqrt{2\pi}}) \leq 1$ instead of $\frac 4 \pi$.
\end{proof}

\subsection{Actions of $\SL_3(\R)$}
The aim of this section is to prove the following theorem.
\begin{thm}\label{thm:FESL3}
Every action by isometries of $\SL_3(\R)$ on a uniformly convex Banach space has a fixed point.
\end{thm}

For $i,j \in \{1,2,3\}$ with $i \neq j$ and $t \in \R$, let $X_{ij}(t)$ denote the elementary matrix with $1$'s on the diagonal, $t$ at the $(i,j)$-entry, and $0$ at the other entries. 

Recall that $E$ denotes a fixed uniformly convex Banach space. Consider an action of $\SL_3(\R)$ by isometries on $E$, and fix $\xi \in E$. Since the elementary subgroups $X_{ij}(\R)$ are exponentially distorted in $\SL_3(\R)$, we know, by Lemma \ref{lem:orbitgrowth}, that there is a constant $M$ (depending on $\xi$) such that for every $i,j$ with $i \neq j$ and $t\in \R$,
\begin{equation}\label{eq:displacement}
\| X_{ij}(t) \cdot \xi-\xi\| \leq M \log(2+|t|).
\end{equation}

For $a_1,a_2,a_3,a_4,a_5,a_6 \geq 0$, define the probability measures
\[
    \mu_{a_1a_2a_3a_4a_5a_6} = X_{12}(\gamma_{a_1})X_{13}(\gamma_{a_2})X_{23}(\gamma_{a_3})X_{21}(\gamma_{a_4})X_{31}(\gamma_{a_5})X_{32}(\gamma_{a_6}),
\]
\[
    \widetilde{\mu}_{a_1a_2a_3a_4a_5a_6} = X_{23}(\gamma_{a_1})X_{13}(\gamma_{a_2})X_{12}(\gamma_{a_3})X_{32}(\gamma_{a_4})X_{31}(\gamma_{a_5})X_{21}(\gamma_{a_6}).
\]
Observe that the measure $\mu_{a_1a_2a_3a_4a_5a_6}$ (resp.~$\widetilde{\mu}_{a_1a_2a_3a_4a_5a_6}$) is obtained by taking clockwise (resp.~counterclockwise) the convolution product of Gaussian measures on the root groups of the root system $A_2$ (which is the root system of $\mathfrak{sl}_3$); see Figure~\ref{picture:A2}.
\\

\begin{figure}
  \center
\begin{tikzpicture}
\draw[stealth-stealth] (180:2cm) node[left] {$X_{21}$}--(0:2cm) node[right] {$X_{12}$};
\draw[stealth-stealth] (240:2cm) node[below left] {$X_{23}$}--(60:2cm)  node[above right] {$X_{32}$};
\draw[stealth-stealth] (300:2cm) node[below right] {$X_{13}$}--(120:2cm)  node[above left] {$X_{31}$};
 \draw[-{Stealth[scale=2]}] (0:3cm) arc (0:-300:3cm);
\draw[-{Stealth[scale=2]}] (-120:3.7cm) arc (-120:180:3.7cm);
\draw (-90:2.7cm) node {$\mu$};
\draw (-45:4cm) node {$\widetilde{\mu}$};
\end{tikzpicture}
\caption{The root system $A_2$ and the measures $\mu_\lambda$ and $\widetilde{\mu}_\lambda$.} \label{picture:A2}
\end{figure}

{\bf Key observation}: All subwords of length $3$ of the measures $\mu_{a_1a_2a_3a_4a_5a_6}$ and $\widetilde{\mu}_{a_1a_2a_3a_4a_5a_6}$, i.e.~the measures
\begin{align*}
    X_{12}(\gamma_{a_1})X_{13}(\gamma_{a_2})X_{23}(\gamma_{a_3}), \quad &X_{13}(\gamma_{a_2})X_{23}(\gamma_{a_3})X_{21}(\gamma_{a_4}), \\
    X_{23}(\gamma_{a_3})X_{21}(\gamma_{a_4})X_{31}(\gamma_{a_5}), \quad 
    &X_{21}(\gamma_{a_4})X_{31}(\gamma_{a_5})X_{32}(\gamma_{a_6}), \\
    X_{23}(\gamma_{a_1})X_{13}(\gamma_{a_2})X_{12}(\gamma_{a_3}), \quad &X_{13}(\gamma_{a_2})X_{12}(\gamma_{a_3})X_{32}(\gamma_{a_4}), \\
    X_{12}(\gamma_{a_3})X_{32}(\gamma_{a_4})X_{31}(\gamma_{a_5}), \quad
    &X_{32}(\gamma_{a_4})X_{31}(\gamma_{a_5})X_{21}(\gamma_{a_6}),
\end{align*}    
are of the form $\psi(\nu_{xyz})$ (with $\nu_{xyz}$ as in subsection~\ref{subsec:heisenberg}) for a continuous homomorphism $\psi \colon H_3(\R) \to \SL_3(\R)$. For example, $X_{12}(\gamma_{a_1})X_{13}(\gamma_{a_2})X_{23}(\gamma_{a_3})$ is the image of $\nu_{a_1a_2a_3}$ under the homomorphism defined by $X(t)\mapsto X_{12}(t)$, $Y(t)\mapsto X_{23}(t)$, and $Z(t) \mapsto X_{13}(t)$. Similarly, $X_{13}(\gamma_{a_2})X_{23}(\gamma_{a_3})X_{21}(\gamma_{a_4})$ is the image of $\nu_{a_2a_3a_4}$ under the homomorphism $X(t)\mapsto X_{13}(t)$, $Y(t)\mapsto X_{21}(t)$ and $Z(t) \mapsto X_{23}(-t)$ (here we use that the Gaussian measure $\gamma_{a_3}$ is symmetric).

For every $\lambda=a_1a_2a_3a_4a_5a_6$, with $a_1,a_2,a_3,a_4,a_5,a_6 \geq 0$, we write
\[
    \xi_\lambda  = \mu_\lambda \cdot \xi \qquad  \textrm{and} \qquad \widetilde{\xi}_\lambda  = \widetilde{\mu}_\lambda \cdot \xi.
\]
For every real number $L>1$, we consider the set of parameters
\[
    P_L = \cup_{s>0} [s,Ls]^6,
\]
i.e.~the set of $6$-tuples of positive numbers, all of whose pairwise ratios are bounded above by $L$.

\begin{lemma} \label{lem:flip_SL3}
    For every $L$, there is a constant $C=C(L) > 0$ such that for every $\lambda = a_1a_2a_3a_4a_5a_6 \in P_L$,
    \[
        \|\xi_{a_1a_2a_3a_4a_5a_6} - \widetilde{\xi}_{a_3a_2a_1a_6a_5a_4}\| \leq C M (a_1 + 1) e^{\max(a_1+a_3-a_2,a_4+a_6-a_5)}.
    \]
\end{lemma}
\begin{proof} 
This follows by applying Lemma~\ref{lem:TV_of_gaussians_H3_affine} and \eqref{eq:displacement} twice. Indeed, we can write $\xi_{a_1a_2a_3a_4a_5a_6} = \nu_1\nu_2\cdot \xi$ and $\widetilde{\xi}_{a_3a_2a_1a_6a_5a_4} = \widetilde{\nu}_1\widetilde{\nu}_2\cdot \xi$, with 
\begin{align*}
    \nu_1 = X_{12}(\gamma_{a_1})X_{13}(\gamma_{a_2})X_{23}(\gamma_{a_3}), \qquad &\nu_2 = X_{21}(\gamma_{a_4})X_{31}(\gamma_{a_5})X_{32}(\gamma_{a_6}), \\
    \widetilde{\nu}_1 = X_{23}(\gamma_{a_3})X_{13}(\gamma_{a_2})X_{12}(\gamma_{a_1}), \qquad &\widetilde{\nu}_2 = X_{32}(\gamma_{a_6})X_{31}(\gamma_{a_5})X_{21}(\gamma_{a_4}).
\end{align*}
By the triangle inequality,
\[
    \|\xi_{a_1a_2a_3a_4a_5a_6} - \widetilde{\xi}_{a_3a_2a_1a_6a_5a_4}\|\leq \| \nu_1 \nu_2 \cdot \xi - \nu_1 \widetilde{\nu}_2 \cdot \xi\| + \| \nu_1 \widetilde{\nu}_2 \cdot \xi - \widetilde{\nu}_1 \widetilde{\nu}_2 \cdot \xi\|.
\]
Using that $\eta\mapsto \nu_1\eta$ is contractive (since the action is by isometries), the first term of the right-hand side is bounded above by 
\begin{align*}
     \|\nu_2\cdot \xi - \widetilde{\nu}_2\cdot \xi\| & \lesssim e^{a_4+a_6-a_5} \max_{|t| \leq e^{a_5}} \|X_{31}(t)\cdot \xi - \xi\|\\
     & \lesssim M (a_5 + 1) e^{a_4+a_6-a_5},
\end{align*}
where the first inequality is Lemma~\ref{lem:TV_of_gaussians_H3_affine}, and the second one is \eqref{eq:displacement}.

By a direct application of Lemma~\ref{lem:TV_of_gaussians_H3_affine}, the second term is bounded above by
\begin{equation} \label{eq:elementarysubgroupandmeasure}
    C e^{a_1+a_3-a_2} \max_{|t|\leq e^{a_2}} \| X_{13}(t) \cdot(\widetilde{\nu}_2 \cdot \xi) - (\widetilde{\nu}_2 \cdot \xi)\|.
\end{equation}
By the bound
\[
    \| X_{13}(t) \cdot(\widetilde{\nu}_2 \cdot \xi) - (\widetilde{\nu}_2 \cdot \xi)\| \leq \| X_{13}(t) \cdot \xi - \xi\| + 2 \|\widetilde{\nu}_2 \cdot \xi-\xi\|,
\]
we deduce that \eqref{eq:elementarysubgroupandmeasure} is  
\[
    \lesssim M (1+a_2+a_4+a_5+a_6).
\]
The lemma now follows from our assumption that $\lambda \in P_L$, which implies $a_5+(a_2+a_4+a_5+a_6) \leq 5L a_1$.
\end{proof}

Let $q=q(E) \in (0,1)$ be given by Proposition~\ref{prop:affine_actions_of_H3R}.

\begin{lemma}\label{lem:local_moves_SL3}
For every $L>1$, there is a constant $C=C(L)>0$ such that if $\varphi$ denotes the function defined by $\func{t} = q^{\sqrt t}$ for $t>0$ and $\infty$ otherwise, then for every $\lambda=a_1a_2a_3a_4a_5a_6$ and $\lambda'$ in $P_L$,
    \begin{equation}\label{eq:move1_SL3} \|\xi_{\lambda} - \xi_{\lambda'}\| \leq C M (1+a_1^2) \func{a_1+a_3-\max(a_2,a_2')} \textrm{ if }\lambda'=a_1a_2'a_3a_4a_5a_6,
    \end{equation} 
    \begin{equation}\label{eq:move2_SL3} \|\xi_{\lambda} - \xi_{\lambda'}\| \leq C M (1+a_1^2) \func{a_2+a_4-\max(a_3,a_3')}\textrm{ if }\lambda'=a_1a_2a_3'a_4a_5a_6,\end{equation}
    \begin{equation}\label{eq:move3_SL3} \|\xi_{\lambda} - \xi_{\lambda'}\| \leq C M (1+a_1^2) \func{a_3+a_5-\max(a_4,a_4')}\textrm{ if }\lambda'=a_1a_2a_3a_4'a_5a_6,\end{equation}
    \begin{equation}\label{eq:move4_SL3} \|\xi_{\lambda} - \xi_{\lambda'}\| \leq C M (1+a_1^2) \func{a_4+a_6-\max(a_5,a_5')}\textrm{ if }\lambda'=a_1a_2a_3a_4a_5'a_6.\end{equation}
\end{lemma}
\begin{proof} As for Lemma~\ref{lem:flip_SL3}, by the key observation, this is just the combination of Proposition~\ref{prop:affine_actions_of_H3R} with the exponential distortion estimate \eqref{eq:displacement}. Observe that here, the upper bound is simpler than in Proposition~\ref{prop:affine_actions_of_H3R} because of our assumption that $\lambda,\lambda' \in P_L$. For example, for the first inequality, Proposition~\ref{prop:affine_actions_of_H3R} would give an upper bound of $\lesssim M q^{\sqrt{a_1+a_3-\max(a_2,a_2')}}(1+|a_2-a_2'|)(1+a_1+a_3)$, but this is less than $Mq^{\sqrt{a_1+a_3-\max(a_2,a_2')}}(1+La_1)(1+(L+1)a_1)\lesssim_L M (1+a_1^2) q^{\sqrt{a_1+a_3-\max(a_2,a_2')}}$ on $P_L$. 
\end{proof}
By symmetry, we obtain the following lemma.
\begin{lemma}\label{lem:local_moves_SL3bis}
With the same notation, we obtain that for every $\lambda=a_1a_2a_3a_4a_5a_6$ and $\lambda'$ in $P_L$,
    \begin{equation}\label{eq:move5_SL3} \|\widetilde{\xi}_{\lambda} - \widetilde{\xi}_{\lambda'}\| \leq C M (1+a_1^2) \func{a_1+a_3-\max(a_2,a_2')} \textrm{ if }\lambda'=a_1a_2'a_3a_4a_5a_6,
    \end{equation} 
    \begin{equation}\label{eq:move6_SL3} \|\widetilde{\xi}_{\lambda} - \widetilde{\xi}_{\lambda'}\| \leq C M (1+a_1^2) \func{a_2+a_4-\max(a_3,a_3')}\textrm{ if }\lambda'=a_1a_2a_3'a_4a_5a_6,\end{equation}
    \begin{equation}\label{eq:move7_SL3} \|\widetilde{\xi}_{\lambda} - \widetilde{\xi}_{\lambda'}\| \leq C M (1+a_1^2) \func{a_3+a_5-\max(a_4,a_4')}\textrm{ if }\lambda'=a_1a_2a_3a_4'a_5a_6,\end{equation}
    \begin{equation}\label{eq:move8_SL3} \|\widetilde{\xi}_{\lambda} - \widetilde{\xi}_{\lambda'}\| \leq C M (1+a_1^2) \func{a_4+a_6-\max(a_5,a_5')}\textrm{ if }\lambda'=a_1a_2a_3a_4a_5'a_6.\end{equation}
\end{lemma}
\begin{proof}
If $\sigma$ is a permutation of $\{1,2,3\}$ and $P_\sigma$ is its permutation matrix, then the automorphism $\alpha_\sigma$ of $\SL_3(\R)$ given by
    \[
        \alpha_\sigma(g) = P_\sigma (g^T)^{-1} P_\sigma^{-1}
    \]
    maps $X_{ij}(t)$ to $X_{\sigma(j)\sigma(i)}(-t)$. Therefore, if $\sigma$ is the permutation $(1\ 3)$, it maps the measure $\mu_{a_1a_2a_3a_4a_5a_6}$ to $\widetilde{\mu}_{a_1a_2a_3a_4a_5a_6}$. This reduces the proof of the lemma to Lemma~\ref{lem:local_moves_SL3}.
\end{proof}

In the light of the above lemmas, we now introduce a natural graph structure on $(0,\infty)^6$.
\begin{defn} \label{defn:graph_SL3}
Let $\mathcal{G}^0=(V^0,A^0)$ be the graph with vertex set $V^0=(0,\infty)^6$ and an edge between $\lambda=a_1a_2a_3a_4a_5a_6$ and $\lambda'$ if one of the following four conditions holds:
\begin{itemize}
    \item $\lambda'=a_1a_2'a_3a_4a_5a_6$ with $\max(a_2,a_2')<a_1+a_3$,
    \item $\lambda'=a_1a_2a_3'a_4a_5a_6$ with $\max(a_3,a_3')<a_2+a_4$,
    \item $\lambda'=a_1a_2a_3a_4'a_5a_6$ with $\max(a_4,a_4')<a_3+a_5$,
    \item $\lambda'=a_1a_2a_3a_4a_5'a_6$ with $\max(a_5,a_5')<a_4+a_6$.
\end{itemize}
\end{defn}
For $\varepsilon \in (0,1)$, we also define the subgraph $\mathcal{G}^0_\varepsilon = (V^0,A^0_\varepsilon)$ of $\mathcal{G}^0=(V^0,A^0)$ with the same vertex set, but with fewer edges, namely only the edges between $\lambda$ and $\lambda'$ satisfying one of the following stronger conditions:
\begin{itemize}
    \item $\lambda'=a_1a_2'a_3a_4a_5a_6$ with $\max(a_2,a_2')\leq(1-\varepsilon)(a_1+a_3)$,
    \item $\lambda'=a_1a_2a_3'a_4a_5a_6$ with $\max(a_3,a_3')\leq(1-\varepsilon)(a_2+a_4)$,
    \item $\lambda'=a_1a_2a_3a_4'a_5a_6$ with $\max(a_4,a_4')\leq(1-\varepsilon)(a_3+a_5)$,
    \item $\lambda'=a_1a_2a_3a_4a_5'a_6$ with $\max(a_5,a_5')\leq(1-\varepsilon)(a_4+a_6)$.
\end{itemize}
Clearly, $E^0=\cup_{\varepsilon>0} E^0_\varepsilon$.

\begin{prop} \label{prop:connectedness}
  Two vertices $\lambda=a_1a_2a_3a_4a_5a_6$ and $\lambda'=a_1'a_2'a_3'a_4'a_5'a_6'$ are in the same connected component of $\mathcal{G}^0$ if and only if $a_1=a_1'$ and $a_6=a_6'$.
\end{prop}
\begin{proof}
    It is clear that $a_1=a_1'$ and $a_6=a_6'$ if $\lambda,\lambda'$ are in the same connected component. The converse follows directly from the following two claims:
        \begin{itemize}
        \item Every $\lambda =a_1a_2a_3a_4a_5a_6$ is in the same connected component as $a_1 m m m m a_6$, where $m=\max(a_1,a_2,a_3,a_4,a_5,a_6)$.
        \item For every $a_1,a_6,m,m' \in (0,\infty)$ satisfying $m<m'<m+a_1$, $a_1 mmmm a_6$ is in the same connected component as $a_1 m'm'm'm' a_6$.
    \end{itemize}

    To prove the first claim, consider first the case $m=a_1$. Then we have the following path:
    \[
        \lambda = ma_2a_3a_4a_5a_6 \sim mma_3a_4a_5a_6 \sim mmma_4a_5a_6 \sim mmmma_5a_6 \sim mmmmma_6.
    \]
    Similarly, if for example $m=a_3$,
    \[
        \lambda = a_1a_2ma_4a_5a_6 \sim a_1mma_4a_5a_6 \sim a_1mmma_5a_6 \sim a_1mmmma_6.
    \]
    The other cases ($m=a_2,a_4,a_5,a_6$) are treated in the same way. Note that we have only used edges of the form $\lambda \sim \lambda'$ if $\lambda'=a_1a_2'a_3a_4a_5a_6$ with $\max(a_2,a_2') < \max(a_1,a_3)$ etc.
    
    Let us consider the second claim. If $m' \leq \max(a_1,a_6)$ then the second claim is a consequence of the first, so we can assume $m' > \max(a_1,a_6)$. Then we have an edge $a_1 mmmm a_6 \sim a_1m' mmma_6$. By the first claim, $a_1m'mmma_6$ is in the same connected component as $a_1m'm'm'm'a_6$. This finishes the proof.
    \end{proof}

Finally, in view of Lemma~\ref{lem:flip_SL3}, it is natural to consider another graph. 
\begin{defn}
Let $\mathcal{G}=(V,A)$ be the graph defined by $V=V^0 \times \{1,2\}$ and an edge beween $(\lambda,i)$ and $(\lambda',j)$ if one of the following conditions holds:
  \begin{itemize}
  \item $i=j$ and $(\lambda,\lambda') \in A^0$,
  \item $i\neq j$ and $\lambda=a_1a_2a_3a_4a_5a_6$ and $\lambda'=a_3a_2a_1a_6a_5a_4$ with $a_1+a_3<a_2$ and $a_4+a_6<a_5$.
  \end{itemize}
Replacing the first condition by $(\lambda,\lambda') \in A^0_\varepsilon$ and the second condition by $a_1+a_3 \leq (1-\varepsilon)a_2$ and $a_4+a_6\leq(1-\varepsilon)a_5$, we obtain a graph that we denote by $\mathcal{G}_\varepsilon$.
\end{defn}

\begin{prop} \label{prop:connectednessbis}
  The graph $\mathcal{G}$ is connected.
\end{prop}
\begin{proof}
    It suffices to show that if $i \neq j$ and $\lambda=a_1a_2a_3a_4a_5a_6,\lambda'=a_1'a_2'a_3'a_4'a_5'a_6' \in (0,\infty)^6$, then $(\lambda,i)$ and $(\lambda',j)$ are in the same connected component. This follows by two uses of Proposition~\ref{prop:connectedness}:
\begin{align*}
    (\lambda,i) &\sim (a_1 \; (2a_1+2a_1') \; a_1' \; a_6' \; (2a_6+2a_6') \; a_6,i) \\ 
        &\sim  (a_1' \; (2a_1+2a_1') \; a_1 \; a_6 \; (2a_6+2a_6') \; a_6',j) \sim (\lambda',j).\qedhere
\end{align*}
\end{proof}
In fact, we will not need the full conclusion of Proposition~\ref{prop:connectednessbis}, only the particular case that $(111111,1)$ and $(222222,1)$ are in the same connected component. Following the proof of the proposition, we can verify that explicitely; see Figure~\ref{figure:example}.
\begin{figure}[h]
\label{figure:example}
\centering
\begin{tikzpicture}
  \newcommand{\smallhistogram}[7]{% #1-#6: heights, #7: position (e.g., (x,y))
    \begin{scope}[shift={#7}]
    \pgfmathsetmacro{\barwidth}{0.15} % Width of each bar
    \pgfmathsetmacro{\totalwidth}{6 * \barwidth} % Total width of the histogram
    \pgfmathsetmacro{\offset}{-\totalwidth/2 + \barwidth/2} % Center the histogram

    \foreach \i/\height in {1/#1, 2/#2, 3/#3, 4/#4, 5/#5, 6/#6} {
      \pgfmathsetmacro{\xpos}{\offset + (\i-1)*\barwidth} % Compute x-position explicitly
      \fill[blue!70] (\xpos, 0) rectangle ++(\barwidth, \height/5);
      % Add height label below the bar
      \node[below, font=\tiny] at (\xpos + \barwidth/2, -0.1) {\pgfmathprintnumber{\height}};
    }
  \end{scope}
}

  \pgfmathsetmacro{\spacement}{2}
  \pgfmathsetmacro{\heighttwo}{-3}
  \pgfmathsetmacro{\heightthree}{-6}
  \pgfmathsetmacro{\deltax}{2}
  \pgfmathsetmacro{\deltay}{3}
  \pgfmathsetmacro{\xpos}{0}
  \pgfmathsetmacro{\ypos}{0}

  \newcommand{\deplacegauche}[1]{
    \pgfmathsetmacro{\xstart}{\xpos-0.5}
    \pgfmathsetmacro{\xend}{\xpos-\deltax+0.6}
    \draw[->] (\xstart,\ypos+0.5) -- node[above,midway] {#1} (\xend,\ypos+0.5);
    \pgfmathsetmacro{\xpos}{\xpos-\deltax}
    }

    \newcommand{\deplacedroite}[1]{
    \pgfmathsetmacro{\xstart}{\xpos+0.6}
    \pgfmathsetmacro{\xend}{\xpos+\deltax-0.5}
    \draw[->] (\xstart,\ypos+0.5) -- node[above,midway] {#1} (\xend,\ypos+0.5);
    \pgfmathsetmacro{\xpos}{\xpos+\deltax}
    }
    \newcommand{\deplacebas}[1]{
    \pgfmathsetmacro{\ystart}{\ypos-0.5}
    \pgfmathsetmacro{\yend}{\ypos-\deltay+1.6}
    \draw[->] (\xpos,\ystart) -- node[left,midway] {#1} (\xpos,\yend);
    \pgfmathsetmacro{\ypos}{\ypos-\deltay}
    }
    \newcommand{\deplacehaut}[1]{
    \pgfmathsetmacro{\ystart}{\ypos+1.6}
    \pgfmathsetmacro{\yend}{\ypos+\deltay-0.5}
    \draw[->] (\xpos,\ystart) -- node[left,midway] {#1} (\xpos,\yend);
    \pgfmathsetmacro{\ypos}{\ypos+\deltay}
    }

  % Draw the 16 histograms in a row manually
  \smallhistogram{3}{3}{3}{3}{3}{3}{(\xpos,\ypos)}
    \draw[dashed] (\xpos+\deltax/2, 3*\deltay/2) -- ((\xpos+\deltax/2, -1.2*\deltay);
  \deplacegauche{2+}
  \smallhistogram{3}{5}{3}{3}{3}{3}{(\xpos,\ypos)}
  \deplacegauche{3+}
  \smallhistogram{3}{5}{7}{3}{3}{3}{(\xpos,\ypos)}
  \deplacegauche{4+}
  \smallhistogram{3}{5}{7}{7}{3}{3}{(\xpos,\ypos)}
  \deplacehaut{5+}
  \smallhistogram{3}{5}{7}{7}{8}{3}{(\xpos,\ypos)}
  \deplacedroite{2+}
  \smallhistogram{3}{8}{7}{7}{8}{3}{(\xpos,\ypos)}
  \deplacedroite{3-}
  \smallhistogram{3}{8}{1}{7}{8}{3}{(\xpos,\ypos)}
  \deplacedroite{4-}
  \smallhistogram{3}{8}{1}{1}{8}{3}{(\xpos,\ypos)}
%  \draw[dashed] (\xpos+\deltax/2, \deltay/2) -- ((\xpos+\deltax/2, -3.2*\deltay);
  \deplacedroite{flip}
  \smallhistogram{1}{8}{3}{3}{8}{1}{(\xpos,\ypos)}
  \deplacebas{3+}
  \smallhistogram{1}{8}{6}{3}{8}{1}{(\xpos,\ypos)}
  \deplacebas{4+}
  \smallhistogram{1}{8}{6}{6}{8}{1}{(\xpos,\ypos)}

  \deplacegauche{flip}
  \smallhistogram{6}{8}{1}{1}{8}{6}{(\xpos,\ypos)}
  \deplacegauche{3+}
  \smallhistogram{6}{8}{6}{1}{8}{6}{(\xpos,\ypos)}
  \deplacegauche{4+}
  \smallhistogram{6}{8}{6}{6}{8}{6}{(\xpos,\ypos)}
  \deplacegauche{2-}
  \smallhistogram{6}{6}{6}{6}{8}{6}{(\xpos,\ypos)}
  \deplacegauche{5-}
  \smallhistogram{6}{6}{6}{6}{6}{6}{(\xpos,\ypos)}
\end{tikzpicture}

\caption{A path of length $15$ between $(333333,1)$ and $(666666,1)$ in $\mathcal{G}$. The dashed line represents the boundary between $V^0\times \{1\}$ on the left and $V^0 \times \{2\}$ on the right. A vertex $(\lambda,i) \in V$ is represented by a histogram with heights $\lambda$. An edge coming from $\mathcal{G}^0$ is oriented and labeled by the corresponding position and a sign, depending on whether the coordinate is increased or decreased. The other edges are labeled \emph{flip}.}
\end{figure}
Theorem \ref{thm:FESL3} will follow from the following important proposition.

\begin{prop} \label{prop:mainSL3}
There is a vector $\xi_\infty$ such that
\[
    \lim_{n \to \infty} \|\xi_{2^n2^n2^n2^n2^n2^n}-\xi_\infty\|=0.
\]
\end{prop}

\begin{proof}[Proof of Theorem~\ref{thm:FESL3} assuming Proposition~\ref{prop:mainSL3}]
Since the vectors $\xi_{aaaaaa}$ satisfy $\lim_{a\to \infty} \|X_{12}(t) \cdot\xi_{aaaaaa} - \xi_{aaaaaa}\|=0$ by Lemma~\ref{lem:TV_of_gaussiansbis} and \eqref{eq:displacement}, $\xi_\infty$ is fixed by $X_{12}(\R)$. As a general consequence of the Mautner phenomenon, we know a priori that the action $\SL_3(\R)\curvearrowright E$ is either proper or it has bounded orbits (see \cite[Theorem 1.4]{MR2541757}). Hence, the action $\SL_3(\R)\curvearrowright E$ has bounded orbits and therefore a fixed point by \cite[Lemma 2.14]{MR2316269}.
\end{proof}

It remains to prove Proposition~\ref{prop:mainSL3}.

\begin{proof}[Proof of Proposition~\ref{prop:mainSL3}]
By Proposition~\ref{prop:connectednessbis} (see also Figure~\ref{figure:example}), the points $(111111,1)$ and $(222222,1)$ are connected by a path in $\mathcal{G}$. Let $k$ denote the length of this path, and let $\varepsilon>0$ and $L>1$ be such that all vertices in this path are in $P_L\times\{1,2\}$ and all edges are in $\mathcal{G}_\varepsilon$. By Example~\ref{figure:example}, we can take $k=15$, $\varepsilon=\frac 1 8$ and $L=8$. By considering the homotheties $(\lambda,i)\mapsto (2^n \lambda,i)$ which preserve $\mathcal{G}$ and $\mathcal{G}_\varepsilon$, we deduce that for every $n$, $(2^n2^n2^n2^n2^n2^n,1)$ and $(2^{n+1}2^{n+1}2^{n+1}2^{n+1}2^{n+1}2^{n+1},1)$ are connected in the graph $\mathcal{G}_\varepsilon$ by a path of length $\leq k$ to all vertices in $P_L\times\{1,2\}$.

Hence, it follows from Lemmas~\ref{lem:flip_SL3},~\ref{lem:local_moves_SL3}, and~\ref{lem:local_moves_SL3bis} that there is $q_1<1$ and $C>1$ such that for every $n \in \Z$, 
\[
    \|\xi_{2^{n}2^{n}2^{n}2^{n}2^{n}2^{n}} - \xi_{2^{n+1}2^{n+1}2^{n+1}2^{n+1}2^{n+1}2^{n+1}}\| \leq C M (1+2^{2n}) q_1^{2^{n/2}},
\]
We immediately deduce that the sequence $(\xi_{2^{n}2^{n}2^{n}2^{n}2^{n}2^{n}})_{n \in \N}$ is Cauchy (and therefore converges). This proves the proposition.
\end{proof}
\begin{rem}
Using the full conclusion of Proposition~\ref{prop:connectednessbis}, the same proof shows that for every cone $K \subset (0,\infty)^6$ with compact basis (i.e., its intersection with the unit sphere of $\R^6$ is compact), we have 
\[
    \lim_{\lambda \in K,\,\lambda \to \infty} \|\xi_\lambda-\xi_\infty\|=0.
\]
\end{rem}

\section{The group $\Sp_4(\R)$}\label{sec:Sp4}

\subsection{Representations and actions of the group $H$} \label{subsection:H}
In this section, we again fix a uniformly convex Banach space $E$.

Consider the $4$-dimensional Lie algebra $\mathfrak{h}$ with basis $\mathfrak{X},\mathfrak{Y},\mathfrak{W},\mathfrak{Z}$ and relations
\[
    [\mathfrak{X},\mathfrak{Y}]=\mathfrak{W}, \quad [\mathfrak{X},\mathfrak{W}]=2\mathfrak{Z}, \quad [\mathfrak{X},\mathfrak{Z}]=[\mathfrak{Y},\mathfrak{Z}]=[\mathfrak{Y},\mathfrak{W}]=0.
\]
\begin{rem}\label{rem:h_contains_h3}
  The Lie algebra $\mathfrak{h}$ contains a subalgebra isomorphic to $\mathfrak{h}_3(\R)$, namely the one spanned by $\mathfrak{X},\mathfrak{W},2\mathfrak{Z}$.
\end{rem}
Let $H$ be the simply connected Lie group with Lie algebra $\mathfrak{h}$, and let $\exp \colon \mathfrak{h} \to H$ denote the corresponding exponential map. For $r,s,t,u \in \R$, the elements $X(u)=\exp(u\mathfrak{X})$, $Y(r) = \exp(r\mathfrak{Y})$, $W(s)=\exp(s\mathfrak{W})$ and $Z(t)=\exp(t\mathfrak{Z})$ satisfy the relations
\[
	[X(u),Y(r)] = Z(-u^2r) W(ur), \qquad [X(u),W(s)]=Z(2us),
\]
and
\[
	[X(u),Z(t)] = [Y(r),W(s)] = [Y(r),Z(t)] = [W(s),Z(t)] = 1.
\]
The group $H$ is a three-step nilpotent group with center $Z(\R)$.

For real numbers $a,b,c,d$, we define the following probability measures on $H$:
\begin{align*}
    \nu_{abcd} := Y(\gamma_a)W(\gamma_b)Z(\gamma_c)X(\gamma_d), \\
    \widetilde{\nu}_{abcd} := X(\gamma_a)Z(\gamma_b)W(\gamma_c)Y(\gamma_d).
\end{align*}

We first introduce some notation.
\begin{defn}
    If $H \curvearrowright E$ is an isometric action, $\xi \in E$ and $L$ is either $X,Y,W,Z$, we write
    \[
        \delta_{L,a}(\xi) = \max_{|t|\leq e^a} \|L(t) \cdot \xi - \xi\|.
    \]
    If $a,b,c,d\in \R$, we write
    \[
        \delta_{abcd}(\xi) = \delta_{Y,a}(\xi) + \delta_{W,b}(\xi) + \delta_{Z,c}(\xi) + \delta_{X,d}(\xi).
    \]
\end{defn}

The main result of this section will compare averages along orbits of an isometric action on $E$ with respect to the measures $\nu_{abcd}$ (and also $\widetilde{\nu}_{dcba}$) when $b$ and $c$ vary. We state it in two distinct results, first for $c$ and then for $b$.
\begin{prop}\label{prop:changec_H}
There is $q < 1$ and $C > 0$ such that for all isometric actions $H \curvearrowright E$, all $a,b,c,c',d \in \R$ with $\max(c,c') \leq b+d$, and all $\xi \in E$,
\[
    \max(\| \nu_{abcd} \cdot \xi - \nu_{abc'd} \cdot \xi\|,\| \widetilde{\nu}_{dcba} \cdot \xi - \widetilde{\nu}_{dc'ba} \cdot \xi\|) \leq C q^{\Delta} (1+|c-c'|) \, \delta_{abcd}(\xi),
\]
where
\[
    \Delta=\sqrt{b+d-\max(c,c')}.
\]
\end{prop}
\begin{proof}
This is an immediate consequence of  Proposition~\ref{prop:affine_actions_of_H3R}. Indeed, by Remark~\ref{rem:h_contains_h3} the elements $X_0(t)=X(t)$, $Y_0(s)=W(s)$ and $Z_0(r)=Z(2r)$ generate a subgroup of $H$ isomorphic to $H_3(\R)$. In this way, we obtain an action by isometries of $H_3(\R)$ on $E$. By the triangle inequality,
\[ \| \nu_{abcd} \cdot \xi - \nu_{abc'd} \cdot \xi\| \leq \| W(\gamma_b) Z(\gamma_c) X(\gamma_d) \cdot \xi -W(\gamma_b) Z(\gamma_{c'}) X(\gamma_d) \cdot \xi\|.\]
The image of the Gaussian measure $\gamma_c$ by the map $r\mapsto r/2$ is the Gaussian measure $\gamma_{c - \log 2}$. In other words,  $Z(\gamma_c) = Z_0(\gamma_{c-\log 2})$. So if we write $c_0=c-\log 2$, $c'_0=c'-\log 2$, the preceding becomes
\[ \| \nu_{abcd} \cdot \xi - \nu_{abc'd} \cdot \xi\| \leq \| Y_0(\gamma_b) Z_0(\gamma_{c_0}) X_0(\gamma_d) \cdot \xi -Y_0(\gamma_b) Z_0(\gamma_{c_0'}) X_0(\gamma_d) \cdot \xi\|.\]
The assumption $\max(c,c')\leq b+d$ implies in particular $\max(c_0,c'_0)\leq b+d$ and by Proposition~\ref{prop:affine_actions_of_H3R}, we know that
\[ \| \nu_{abcd} \cdot \xi - \nu_{abc'd} \cdot \xi\| \leq C q^\Delta (1+|c-c'|) \max_{r \leq e^{b+d}} \|Z_0(2r)\cdot \xi - \xi\|.\]
The other inequality is proved the same way: If $\eta = Y(\gamma_a) \xi$, then
\[ \| \widetilde{\nu}_{dcba} \cdot \xi - \widetilde{\nu}_{dc'ba} \cdot \xi\|  = \| X_0(\gamma_d) Z_0(\gamma_{c_0}) Y_0(\gamma_b)  X_0(\gamma_d) \cdot \eta - X_0(\gamma_d)Z_0(\gamma_{c'_0}) Y_0(\gamma_b)  \cdot \eta\|.\]
By Proposition~\ref{prop:affine_actions_of_H3R} this is bounded above by
\[ C q^\Delta (1+|c-c'|) \max_{r \leq e^{b+d}} \|Z_0(2r)\cdot \eta - \eta\| \leq C q^\Delta (1+|c-c'|) \max_{r \leq e^{b+d}} \|Z_0(2r)\cdot \xi- \xi\|. \]
The last inequality uses that $Y(\R)$ commutes with $Z_0(\R)$. 
The proposition follows, because
\[ \max_{|r| \leq e^{b+d}} \|Z_0(2r)\cdot \xi - \xi\| = \max_{|s| \leq e^b,|t|\leq e^d} \| [W(s),X(t)] \cdot \xi - \xi\| \leq 2 \delta_{abcd}(\xi).\qedhere\]
\end{proof}
\begin{prop}\label{prop:nunutildechangeb}
There is $q < 1$ and $C > 0$ such that for all isometric actions $H \curvearrowright E$, all $a,b,b',c,d \in \R$ with $\max(b,b') \leq \min(a+d, \frac{a+c}{2})$, and all $\xi \in E$,
\[
    \| \nu_{abcd} \cdot \xi - \nu_{ab'cd} \cdot \xi\| \leq C q^{\Delta} (1+|b-b'|)\left(\delta_{X,d}(\xi) + \delta_{Y,a}(\xi)\right),
\]
and 
\[
    \| \widetilde{\nu}_{dcba} \cdot \xi - \widetilde{\nu}_{dcb'a} \cdot \xi\| \leq C q^{\Delta} (1+|b-b'|)\left(\delta_{X,d}(\xi) + \delta_{Y,a}(\xi)\right),
\]
where
\[
    \Delta=\sqrt{\min(a+d-b,a+d-b',a+c-2b,a+c-2b')}.
\]
\end{prop}
The rest of this subsection is devoted to the proof of Proposition~\ref{prop:nunutildechangeb}. Our first auxiliary lemma is analogous to Lemma~\ref{lem:TV_of_gaussians_H3_affine}; it will also be used in the next subsection.

\begin{lemma}\label{lem:flip_H}
    There exists $C > 0$ such that for all isometric actions $H \curvearrowright E$, all $a,b,c,d \in \R$, and all $\xi \in E$,
	\[
		\|\nu_{abcd}\cdot \xi - \widetilde{\nu}_{dcba}\cdot \xi\| \leqslant C \big( e^{a+d-b} + e^{b+d-c} \big) \, \delta_{abcd}(\xi).
	\]
\end{lemma}
To prove this lemma, we need two sublemmas.
\begin{lemma} \label{lem:TV_H4_1}
In the setting of Lemma~\ref{lem:flip_H},
    \[
        \|X(\gamma_d) W(\gamma_b) Z(\gamma_c)\cdot \xi - Z(\gamma_c) W(\gamma_b) X(\gamma_d) \cdot \xi\| \leq 4 e^{b+d-c} \delta_{Z,c}(\xi).
    \]
\end{lemma}
\begin{proof}
This follows from Lemma~\ref{lem:TV_of_gaussians_H3_affine}. Indeed the group generated by $X(\R),W(\R),Z(\R)$ is isomorphic to $H_3(\R)$ (see Remark~\ref{rem:h_contains_h3}), therefore the action of $H$ induces an action of $H_3(\R)$ and
\begin{equation} \label{eq:xwzzwx}
    \|X(\gamma_d) W(\gamma_b) Z(\gamma_c)\cdot \xi - Z(\gamma_c) W(\gamma_b) X(\gamma_d) \cdot \xi\| = \| \nu_{d(c-\log 2)b}\cdot \xi - \widetilde{\nu}_{b(c-\log 2)d}\cdot \xi\|.
\end{equation}
Note that on the right-hand side, we use the measures of the form $\nu_{xyz}$ and $\widetilde{\nu}_{zyx}$ from Section \ref{sec:SL3} (with three indices). Lemma~\ref{lem:TV_of_gaussians_H3_affine} shows that \eqref{eq:xwzzwx} is bounded above by $2 e^{b+d-(c-\log 2)} \delta_{Z,c-\log 2}(\xi) \leq 4 e^{b+d-\log 2} \delta_{Z,c}(\xi)$.
\end{proof}
\begin{lemma} \label{lem:TV_H4_2}
In the setting of Lemma~\ref{lem:flip_H},
        \[
            \|\widetilde{\nu}_{dcba} \cdot \xi - Y(\gamma_a) X(\gamma_d) W(\gamma_b) Z(\gamma_c)\cdot \xi\| \leq 2e^{a+d-b}\delta_{W,b}(\xi) + 2 e^{a+2d-c}\delta_{Z,c}(\xi).
        \]
\end{lemma}
\begin{proof}
    By the triangle inequality, using that $Y(\R)$ centralizes the group generated by $W(\R)$ and $Z(\R)$, this norm is bounded above by
    \[
        \E_{u\sim \gamma_d,r\sim \gamma_a} \|X(u) Y(r) \cdot (W(\gamma_b) Z(\gamma_c)\cdot \xi) - Y(r) X(u)  \cdot (W(\gamma_b) Z(\gamma_c)\cdot \xi)\|.
    \]
    By the commutation relation $[X(u),Y(r)] = Z(-u^2r)W(ur)$, this is equal to
    \[
        \E_{u\sim \gamma_d,r\sim \gamma_a} \|W(\gamma_{ur,b})Z(\gamma_{-u^2r,c})\cdot \xi - W(\gamma_b) Z(\gamma_c) \cdot \xi \|,
    \]
    which by the triangle inequality and Lemma~\ref{lem:TV_of_gaussiansbis} is bounded above by
    \begin{align*}
        &2\E_{u\sim \gamma_d,r\sim \gamma_a} (|ur|e^{-b}\delta_{W,b}(\xi) + u^2 |r|e^{-c}\delta_{Z,c}(\xi)) \\ &= \frac{4} {\pi} e^{a+d-b} \delta_{W,b}(\xi) + \frac{2\sqrt 2}{\sqrt \pi} e^{a+2d-c}\delta_{Z,c}(\xi).
    \end{align*}
\end{proof}

\begin{proof}[Proof of Lemma~\ref{lem:flip_H}]
We can assume that $a+d-b<0$ and $b+d-c<0$; otherwise the inequality is easy by the triangle inequality and Lemma~\ref{lem:TV_of_gaussiansbis}:
\[ \|\nu_{abcd}\cdot \xi - \widetilde{\nu}_{dcba}\cdot \xi\| \leq \|\nu_{abcd}\cdot \xi - \xi\| + \|\xi-\widetilde{\nu}_{dcba}\cdot \xi\| \leq 2\left(1+\sqrt{\frac{2}{\pi}}\right) \delta_{abcd}(\xi).\]
Convolution by a fixed probability measure decreases the total variation norm, so Lemma~\ref{lem:TV_H4_1} implies
\[
    \| \nu_{abcd}\cdot \xi - Y(\gamma_a) X(\gamma_d) W(\gamma_b) Z(\gamma_c)\cdot \xi\|_{TV} \leq 4 e^{b+d-c}\delta_{Z,c}(\xi).
\]
On the other hand, Lemma~\ref{lem:TV_H4_2} yields 
\[
    \| \widetilde{\nu}_{dcba}\cdot \xi - Y(\gamma_a) X(\gamma_d) W(\gamma_b) Z(\gamma_c)\cdot \xi\| \leq 2e^{a+d-b}\delta_{W,b}(\xi) + 2e^{a+2d-c}\delta_{Z,c}(\xi).
\]
Altogether, we obtain
\[
    \|\nu_{abcd} - \widetilde{\nu}_{dcba}\|_{TV} \leq 2 e^{b+d-c}\delta_{W,b}(\xi) + (4e^{a+d-b} + 2e^{a+2d-c})\delta_{Z,c}(\xi).
\]
The lemma now follows with $C=6$, because
\[
    e^{a+2d-c} = e^{a+d-b} e^{b+d-c} \leq e^{a+d-b}+ e^{b+d-c}
\]
by our assumption $a+d-b<0$ and $b+d-c<0$.
\end{proof}

The next lemma is a consequence of Proposition~\ref{prop:uniform_convexity_iteratedR}.
\begin{lemma} \label{lem:xyestimate}
There exists $q_0 < 1$ such that for all isometric representations \mbox{$\pi \colon H \to \O(E)$}, all $a,d > 0$, and all $\xi \in E$,
\[
    \max( \|\pi(X(\gamma_d)Y(\gamma_a))\xi\| , \|\pi(Y(\gamma_a)X(\gamma_d))\xi\| ) \leq \max(q_0 \|\xi\|,2\|\pi(W(\gamma_{a+d}))\xi\|).
\]
\end{lemma}
%\begin{rem}
%  The proof of Lemma~\ref{lem:xyestimate} will give a slightly better estimate (which we will not use), namely with the left-hand side replaced by
%\[
%    \max(q_0 \|\xi\|,2\|\pi(Z(\gamma_{a+2d})W(\gamma_{a+d}))\xi\|).
%\]
%\end{rem}

The proof of Lemma \ref{lem:xyestimate} relies on a simple commutator estimate, which we record for later use in the following sublemma. 
\begin{lemma} \label{lem:commutator_cocycle_growth}
    If $H\curvearrowright E$ is an isometric action, then for all $a,d \in \R$,
    \begin{align*}
    \delta_{W,a+d}(\xi) &\leq 4 \delta_{X,d}(\xi) + 4\delta_{Y,a}(\xi),\\
    \delta_{Z,a+2d}(\xi) &\leq 4 \delta_{X,d}(\xi) + 4\delta_{Y,a}(\xi).
    \end{align*}
\end{lemma}
\begin{proof}
    Set $C'=\delta_{X,d}(\xi) + \delta_{Y,a}(\xi)$. For all $|u|\leq e^d$ and $|r|\leq e^a$, by the commutation relation $[X(u),Y(r)] = Z(-u^2r) W(ur)$, we deduce
    \begin{align*}
        W(2ur) &= [X(u),Y(r)] \, [X(-u),Y(-r)], \\
        Z(2u^2r) &= [X(u),Y(-r)] \, [X(-u),Y(-r)].
    \end{align*}
    Hence,
    \begin{align*}
    \| W(2ur)\cdot \xi - \xi\| &\leq 4 C'\\
    \| Z(2u^2r)\cdot \xi - \xi\| & \leq 4 C'.
    \end{align*}
The lemma follows by taking the supremum over $u,r$.
\end{proof}
\begin{proof}[Proof of Lemma~\ref{lem:xyestimate}]
  Let $\varepsilon = \frac{1}{16\left(1+\sqrt{\frac{2}{\pi}}\right)}$, and let $\delta>0$ be given by Proposition~\ref{prop:uniform_convexity_iteratedR} for $n=2$ and this $\varepsilon$. We prove the lemma with $q_0=1-\delta$. Let $\pi$, $a$, $d$, and $\xi$ be given. Suppose $\|\xi\|=1$ and
\[
    \max( \|\pi(X(\gamma_d)Y(\gamma_a))\xi\| , \|\pi(Y(\gamma_a)X(\gamma_d))\xi\| ) \geq q_0.
\]
By Proposition~\ref{prop:uniform_convexity_iteratedR}, we have $\|\pi(X(u)) \xi - \xi\| \leq \varepsilon$ for every $|u| \leq e^d$ and $\|\pi(Y(r)) \xi - \xi\| \leq \varepsilon$ for every $|r| \leq e^a$. By Lemma~\ref{lem:commutator_cocycle_growth}, we deduce $\|\pi(W(s))\xi-\xi\| \leq 8\varepsilon$ for all $|s| \leq e^{a+d}$, and therefore, by Lemma~\ref{lem:from_pointwise_to_gaussian}, we have
\[
    \|\pi(W(\gamma_{a+d}))\xi-\xi\| \leq 8\left(1+\sqrt{\frac 2 \pi}\right)\varepsilon = \frac{1}{2}.
\]
Hence, $\|\pi(W(\gamma_{a+d}))\xi\|\geq \frac 1 2$, which proves the lemma.
\end{proof}

\begin{lemma} \label{lem:calculusbis}
For all $a,b,c,d$ and $\Delta$ as in Proposition~\ref{prop:nunutildechangeb} with $\Delta \geq 7$, there exists an integer $n$ with $\frac{\Delta}{3} \leq n \leq \Delta$ and real numbers $a_1,\dots,a_n$, $d_1,\dots,d_n$ such that
\begin{enumerate}[(i)]
    \item\label{item1bis} $a_i+d_i \geq b+\Delta$ for all $i=1,\ldots,n$,
    \item\label{item2bis} $a_i + d_{i+1} \leq b-\Delta$ for all $i=1,\ldots,n-1$,
    \item\label{item3bis} $b+d_{i+1} \leq c-\Delta$ for all $i=1,\ldots,n-1$,
    \item\label{item4bis} $\sum_{i=1}^n e^{2a_i} = e^{2a}$ and $\sum_{i=1}^n e^{2d_i}=e^{2d}$.
\end{enumerate}
\end{lemma}
\begin{proof}
Set $n=\lfloor \frac{\Delta}{3} \rfloor+1$, $a_0=\max(b-d,2b-c)$, $a_i = a_0 + 2i\Delta$ and $d_i = b+\Delta-a_0-2i\Delta$ for $i=1,\ldots,n$. It is a small computation to verify (\ref{item1bis}), (\ref{item2bis}) and (\ref{item3bis}). In a similar way as in Lemma~\ref{lem:calculus}, we obtain (\ref{item4bis}) by increasing $a_n$ and $d_1$.
\end{proof}

\begin{proof}[Proof of Proposition~\ref{prop:nunutildechangeb}]
Consider an arbitrary isometric action $H \curvearrowright E$. In the same way as in Proposition \ref{prop:affine_actions_of_H3R}, we can assume $\Delta \geq 7$ and $b < b' \leq b+1$. Moreover, using that $Z(\R)$ is contained in the center of $H$, we see that the conclusion of the proposition for some $c$ implies its validity for all larger $c$. We can therefore assume that $\frac{a+c}{2} \leq a+d$, that is, $c \leq a+2d$. In that case, we know by Lemma~\ref{lem:commutator_cocycle_growth} that $\delta_{abcd}(\xi) \leq 9(\delta_{X,d}(\xi) + \delta_{Y,a}(\xi))$. With these reductions we are left to prove 
\begin{equation} \label{eq:change_b_affineH1}
    \| \nu_{abcd} \cdot \xi - \nu_{ab'cd} \cdot \xi\| \lesssim q^{\Delta}\delta_{abcd}(\xi),
\end{equation}
and
\begin{equation} \label{eq:change_b_affineH2}
\| \widetilde{\nu}_{dcba} \cdot \xi - \widetilde{\nu}_{dcb'a} \cdot \xi\| \lesssim q^{\Delta}\delta_{abcd}(\xi).
\end{equation}

Let $n$, $a_1,\ldots,a_n$, and $d_1,\ldots,d_n$ be given by Lemma~\ref{lem:calculusbis}. We first prove \eqref{eq:change_b_affineH1}. To this end, for $0 \leq k \leq n$, define $a_k^\prime$ and $d_k^\prime$ by $\exp(2a_k^\prime)=\sum_{j=1}^{k} \exp(2a_j)$ and $\exp(2d_k^\prime)=\sum_{j=1}^{k} \exp(2d_j)$ (so $a_0'=d_0'=-\infty$), and set
\[
    \eta_k=\nu_{a_k^\prime bcd_k^\prime} \cdot \xi - \nu_{a_k^\prime b^\prime cd_k^\prime} \cdot \xi.
\]
By Lemma~\ref{lem:calculusbis}.\eqref{item4bis}, we have $\eta_n = \nu_{abcd} \cdot \xi - \nu_{ab'cd} \cdot \xi$, and $\eta_0=W(\gamma_b)Z(\gamma_c)\cdot\xi-W(\gamma_{b'})Z(\gamma_c)\cdot\xi$. By the convolution properties of Gaussian measures, for $1 \leq k \leq n$, we have
\[
    \nu_{a_k^\prime bcd_k^\prime} \cdot \xi =  Y(\gamma_{a_k})\nu_{a_{k-1}^\prime bcd_{k}} X(\gamma_{d_{k-1}^\prime}) \cdot \xi,
\]
so by Lemma~\ref{lem:flip_H} and Lemma~\ref{lem:calculusbis}.\eqref{item2bis} and \eqref{item3bis},
\begin{align*}
    \nu_{a_k^\prime bcd_k^\prime} \cdot \xi
    &= Y(\gamma_{a_k})\widetilde{\nu}_{d_k cb a_{k-1}^\prime}X(\gamma_{d_{k-1}^\prime}) \cdot \xi + O(\sqrt{\Delta}e^{-\Delta})\delta_{abcd}(\xi) \\
    &= Y(\gamma_{a_k})X(\gamma_{d_k}) \nu_{a_{k-1}^\prime bcd_{k-1}^\prime} \cdot \xi + O(\sqrt{\Delta}e^{-\Delta})\delta_{abcd}(\xi).
\end{align*}
In the first equality, we have used that $e^{a'_{k-1} + d_k -b} \leq  \sqrt{n} \max_{j<k} e^{a_j+d_k-b} \leq \sqrt{\Delta} e^{-\Delta}$.

The same estimate holds for $b$ replaced by $b'$, so
\[
    \eta_k = \pi(Y(\gamma_{a_k})X(\gamma_{d_k}))\eta_{k-1} + O(\sqrt{\Delta}e^{-\Delta})\delta_{abcd}(\xi),
\]
where $\pi$ is the isometric representation underlying the action. Applying Lemma~\ref{lem:xyestimate} and Lemma~\ref{lem:calculusbis}.\eqref{item1bis}, we obtain
\[
    \|\eta_k\| \leq  \max(q_0\|\eta_{k-1}\|,2\|\pi(W(\gamma_{b+\Delta}))\eta_{k-1}\|)+O(\sqrt{\Delta}e^{-\Delta})\delta_{abcd}(\xi).
\]
Since
\begin{align*}
    &\pi(W(\gamma_{b+\Delta}))\eta_{k-1} = W(\gamma_{b+\Delta}) \nu_{a_{k-1}^\prime bcd_{k-1}^\prime} \cdot \xi - W(\gamma_{b+\Delta}) \nu_{a_{k-1}^\prime b^\prime cd_{k-1}^\prime} \cdot \xi \\
&\,= \pi(Y(\gamma_{a_{k-1}}^\prime)Z(\gamma_c)) ( W(\gamma_{b + \Delta})W(\gamma_b)X(\gamma_{d_{k-1}^\prime}) \cdot \xi - W(\gamma_{b + \Delta})W(\gamma_{b^\prime})X(\gamma_{d_{k-1}^\prime})\cdot \xi ),
\end{align*}
we obtain, using Lemma~\ref{lem:TV_of_gaussiansbis}, that
\begin{align*}
    \|\pi(W(&\gamma_{b+\Delta}))\eta_{k-1}\| \\ & \leq \frac{e^{2b'}-e^{2b}}{e^{2b+2\Delta}} \max_{|s|\leq \sqrt{ e^{2b+2\Delta} + e^{2b}}} \|W(s)\cdot (X(\gamma_{d_{k-1}^\prime})\cdot \xi) - X(\gamma_{d_{k-1}^\prime})\cdot \xi\|\\
    &\lesssim \exp(-2\Delta) \, \delta_{abcd}(\xi).
\end{align*}
The second inequality is by Lemma~\ref{lem:commutator_cocycle_growth}. Indeed (recall $a+d \geq b+\Delta^2$ and $\Delta\geq 7$), we have $\sqrt{ e^{2b+2\Delta} + e^{2b}} \leq e^{a+d}$. As a conclusion, we obtain
\[
    \| \nu_{abcd} \cdot \xi - \nu_{ab'cd} \cdot \xi\| = \|\eta_n\| \lesssim \left(q_0^n + \sqrt{\Delta} e^{-\Delta} \right) \delta_{abcd}(\xi),
\]
so \eqref{eq:change_b_affineH1} follows as soon as $q \leq q_0^{\frac{1}{3}}$ and $q<e^{-1}$.

We now prove \eqref{eq:change_b_affineH2}. To this end, for $1 \leq k \leq n+1$, define $\tilde{a}_k$ and $\tilde{d}_k$ by $e^{2\tilde{a}_k}=\sum_{j=k}^n e^{2a_j}$ and $e^{2\tilde{d}_k}=\sum_{j=k}^n e^{2d_j}$ (so $\tilde{a}_{n+1} = \tilde{d}_{n+1}=-\infty$). Note that $\tilde{a}_k \neq a_k'$ and $\tilde{d}_k \neq d'_k$; the inequalities are reversed. Define 
\[
    \eta_k = \widetilde{\nu}_{\tilde{d}_k cb \tilde{a}_k} \cdot \xi - \widetilde{\nu}_{\tilde{d}_k cb' \tilde{a}_k} \cdot \xi,
\]
so that $\eta_1 = \widetilde{\nu}_{dcba} \cdot \xi - \widetilde{\nu}_{dcb'a} \cdot \xi $ and $\eta_{n+1} = Z(\gamma_c) W(\gamma_b)\cdot \xi - Z(\gamma_c) W(\gamma_{b'})\cdot \xi$. For $1 \leq k \leq n$, we have, similar to the first case, that
\begin{align*}
    \widetilde{\nu}_{\tilde{d}_k cb \tilde{a}_k} \cdot \xi  & = X(\gamma_{d_k}) \widetilde{\nu}_{\tilde{d}_{k+1} cb a_k} Y(\gamma_{\tilde{a}_{k+1}}) \cdot \xi\\
& = X(\gamma_{d_k}) \nu_{a_k bc \tilde{d}_{k+1}} Y(\gamma_{\tilde{a}_{k+1}}) \cdot \xi + O(\sqrt{\Delta} e^{-\Delta})\delta_{abcd}(\xi)\\
& = X(\gamma_{d_k}) Y(\gamma_{a_k}) \cdot \widetilde{\nu}_{\tilde{d}_{k+1} cb \tilde{a}_{k+1}} \cdot \xi + O(\sqrt{\Delta} e^{-\Delta})\delta_{abcd}(\xi).
\end{align*}
By the same inequality for $b'$, Lemma~\ref{lem:xyestimate} and Lemma~\ref{lem:calculusbis}.\eqref{item1bis}, we obtain
\[
    \|\eta_k\| \leq  \max(q_0\|\eta_{k+1}\|,2\|\pi(W(\gamma_{b+\Delta}))\eta_{k+1}\|)+O(\sqrt{\Delta}e^{-\Delta})\delta_{abcd}(\xi).
\]
Note that this inequality is also valid for $k=n$ (even without the $O(\cdot)$-summand). Let $\zeta=W(\gamma_b) Y(\gamma_{\tilde{a}_{k+1}}) \cdot \xi - W(\gamma_{b'}) Y(\gamma_{\tilde{a}_{k+1}}) \cdot \xi$, so that $\eta_{k+1} = \pi(X(\gamma_{\tilde{d}_{k+1}}) Z(\gamma_c)) \zeta$. By Lemma \ref{lem:TV_H4_1}, we obtain the following bound:
\begin{align*}
    \|\pi(W(\gamma_{b+\Delta}))\eta_{k+1}\| & = \| \pi(W(\gamma_{b+\Delta}) X(\gamma_{\tilde{d}_{k+1}}) Z(\gamma_c)) \zeta\|\\
    &\leq \|\pi(X(\gamma_{\tilde{d}_{k+1}}) Z(\gamma_c)W(\gamma_{b+\Delta})) \zeta\| + O(\sqrt{\Delta} e^{-\Delta})\|\zeta\| \\
    &\leq \|W(\gamma_{b+\Delta} \gamma_b)\cdot \xi - W(\gamma_{b+\Delta} \gamma_{b'})\cdot\xi\| + O(\sqrt{\Delta} e^{-\Delta})\|\zeta\|\\
    &\lesssim e^{-2\Delta} \max_{|s| \leq \sqrt{e^{2b+2\Delta}+e^{2b}}} \|W(s)\cdot\xi-\xi\|+\sqrt{\Delta} e^{-\Delta}\|\zeta\|,
\end{align*}
where the last inequality follows from Lemma~\ref{lem:TV_of_gaussiansbis}. Using Lemma~\ref{lem:commutator_cocycle_growth} and \eqref{eq:tail_of_gaussians}, we deduce
\[
    \|\pi(W(\gamma_{b+\Delta}))\eta_{k+1}\| \lesssim \sqrt{\Delta} e^{-\Delta} \delta_{abcd}(\xi).
\] 
Putting everything together, we obtain
\[
    \|\eta_k\| \leq  q_0\|\eta_{k+1}\| + O(\sqrt{\Delta}e^{-\Delta})\delta_{abcd}(\xi),
\]
which implies $\|\eta_1\| \lesssim (q_0^n + \sqrt{\Delta}e^{-\Delta})\delta_{abcd}(\xi)$, because
\[
    \|\eta_{n+1}\| \leq \| W(\gamma_b)\cdot \xi - W(\gamma_{b'}) \cdot \xi\| \lesssim \delta_{abcd}(\xi)
\]
by Lemma~\ref{lem:commutator_cocycle_growth}. This concludes the proof of \eqref{eq:change_b_affineH2} and of the proposition.
\end{proof}

\subsection{Actions of $\Sp_4(\R)$}
The aim of this section is to prove the following theorem.
\begin{thm}\label{thm:FESp4}
Every action by isometries of $\Sp_4(\R)$ on a uniformly convex Banach space has a fixed point.
\end{thm}

Consider the Lie group
\[
	\Sp_4(\R) = \{ g \in \mathrm{GL}_4(\R) \mid g^TJg=J \}, \quad \textrm{where } J = \begin{pmatrix} 0 & I_2 \\ -I_2 & 0
\end{pmatrix},
\]
and its Lie algebra
\[
	\mathfrak{sp}_4(\R) = \{ M \in \mathrm{Mat}_4(\R) \mid M^T J + J M = 0\}.
\]
As a Cartan subalgebra of $\mathfrak{sp}_4(\R)$, we choose the subalgebra
\[
	\mathfrak{a} = \left\{ D(a,b)=\begin{pmatrix} a & 0 & 0 & 0\\
                                    0 & b & 0 & 0\\
                                    0 & 0 & -a & 0\\
                                    0 & 0 & 0 & -b
\end{pmatrix} \Bigg\vert \; a,b \in \R\right\}.
\]

The associated root system $C_2$ is given by $\{\pm 2\alpha, \pm 2\beta, \pm \alpha \pm \beta\}$, where $\alpha,\beta \in \mathfrak{a}^*$ are given by $\alpha(D(a,b))=a$ and $\beta(D(a,b))=b$. Each root space has dimension $1$, and we choose the following specific vectors $\mathfrak{X}_\phi$ (where $\phi$ is a root) spanning the associated root space:
\begin{align*}
\mathfrak{X}_{2\beta} = \begin{pmatrix} 0 & 0 & 0 & 0\\
                       0 & 0 & 0 & 1\\
                       0 & 0 & 0 & 0\\
                       0 & 0 & 0 & 0
\end{pmatrix}, \qquad \mathfrak{X}_{\alpha+\beta} &= \begin{pmatrix} 0 & 0 & 0 & 1\\
                       0 & 0 & 1 & 0\\
                       0 & 0 & 0 & 0\\
                       0 & 0 & 0 & 0
\end{pmatrix},\\
\mathfrak{X}_{2\alpha} = \begin{pmatrix} 0 & 0 & 1 & 0\\
                       0 & 0 & 0 & 0\\
                       0 & 0 & 0 & 0\\
                       0 & 0 & 0 & 0
\end{pmatrix}, \qquad \mathfrak{X}_{\alpha - \beta} &= \begin{pmatrix} 0 & 1 & 0 & 0\\
                       0 & 0 & 0 & 0\\
                       0 & 0 & 0 & 0\\
                       0 & 0 & -1 & 0
\end{pmatrix},
\end{align*}
and $\mathfrak{X}_{-\phi} = - \mathfrak{X}_\phi^T$.

For every root $\phi$ and $t \in \R$, let $X_\phi(t) = \exp{(t\mathfrak{X_\phi})} = I_4 + t\mathfrak{X}_\phi$. For the sake of brevity, we index the roots by $\Z/8\Z$ in the clockwise order (see Figure \ref{picture:C2}):
\begin{align*}
\phi_1 &= 2\beta, &\phi_2 &= \beta +\alpha, &\phi_3 &=2\alpha, &\phi_4 &=\alpha - \beta, \\
\phi_5 &= -2\beta, &\phi_6 &=-\beta - \alpha, &\phi_7 &= -2\alpha, &\phi_8 &= \beta - \alpha,
\end{align*}
and we write $X_i(t)$ for $X_{\phi_i}(t)$.

Below we use that for every $i$, the groups $X_i(\R),X_{i+1}(\R),X_{i+2}(\R),X_{i+3}(\R)$ generate a copy of $H$ in $\Sp_4(\R)$. If $i$ is odd, such an identification is given by
\[
    X(t) \mapsto X_{i+3}(t), \; Y(t) \mapsto X_i(t), \; W(t) \mapsto X_{i+1}(t), \; Z(t)\mapsto X_{i+2}(t).
\] 
If $i$ is even, such an identification is given by
\[
    X(t) \mapsto X_{i}(t), \; Y(t) \mapsto X_{i+3}(t), \; W(t) \mapsto X_{i+2}(-t), \; Z(t)\mapsto X_{i+1}(t).
\]
These embeddings are all conjugates of the standard embedding (for $i=1$) by an inner automorphism of $\Sp_4(\R)$.

Consider an action of $\Sp_4(\R)$ by isometries on our fixed uniformly convex Banach space $E$, and fix $\xi \in E$. The root subgroups being exponentially distorted in $\Sp_4(\R)$, we know by Lemma~\ref{lem:orbitgrowth} that there is a constant $M$ (depending on $\xi$) such that for every $\phi \in C_2$ and $t\in \R$,
\begin{equation}\label{eq:displacement_Sp4}
\| X_{\phi}(t) \cdot \xi-\xi\| \leq M \log(2+|t|).
\end{equation}

For $\lambda = (a_1,\dots,a_8) \in [0,\infty]^8$, set $|\lambda| = \max_i a_i$, and define the following probability measures:
\[
    \mu_{\lambda} = X_1(\gamma_{a_1}) X_2(\gamma_{a_2}) X_3(\gamma_{a_3}) X_4(\gamma_{a_4})X_5(\gamma_{a_5})X_6(\gamma_{a_6})X_7(\gamma_{a_7})X_8(\gamma_{a_8}) ,
\]
\[
    \widetilde{\mu}_{\lambda} = X_4(\gamma_{a_1})X_3(\gamma_{a_2}) X_2(\gamma_{a_3})X_1(\gamma_{a_4})X_8(\gamma_{a_5})X_7(\gamma_{a_6})X_6(\gamma_{a_7})X_5(\gamma_{a_8}).
\]
We also define the following elements of $E$:
\[
    \xi_\lambda  = \mu_\lambda \cdot \xi \qquad  \textrm{and} \qquad \widetilde{\xi}_\lambda  = \widetilde{\mu}_\lambda \cdot \xi.
\]

Observe that the measure $\mu_{\lambda}$ (resp.~$\widetilde{\mu}_{\lambda}$) is obtained by taking clockwise (resp.~counterclockwise) the convolution product of Gaussian measures on the root groups of the root system $C_2$; see Figure \ref{picture:C2}.\\
\begin{figure}
  \center
\begin{tikzpicture}
\draw[stealth-stealth] (-1,-1) node[below left] {$X_{6}$}--(1,1) node[above right] {$X_{2}$};
\draw[stealth-stealth] (0,-2) node[below] {$X_{5}$}--(0,2) node[above] {$X_{1}$};
\draw[stealth-stealth] (-1,1) node[above left] {$X_{8}$}--(1,-1) node[below right] {$X_{4}$};
\draw[stealth-stealth] (-2,0) node[left] {$X_{7}$}--(2,0) node[right] {$X_{3}$};
\draw[-{Stealth[scale=2]}] (90:3cm) arc (90:-225:3cm);
\draw[-{Stealth[scale=2]}] (-45:4cm) arc (-45:270:4cm);
\draw (70:2.7cm) node {$\mu$};
\draw (-20:4.5cm) node {$\tilde{\mu}$};
\end{tikzpicture}
\caption{The root system $C_2$ and the measures $\mu_\lambda$ and $\widetilde{\mu}_\lambda$.} \label{picture:C2}
\end{figure}

{\bf Key observation}: All subwords of length $4$ defining the measures $\mu_{\lambda}$ and $\widetilde{\mu}_{\lambda}$ are of the form $\psi(\nu_{a'b'c'd'})$ or $\psi(\widetilde{\nu}_{a'b'c'd'})$ for a continuous homomorphism $\psi \colon H \to \Sp_4(\R)$, depending on whether the first word is $X_i(\cdot)$ for an odd or even $i$. 

The main result of this section, i.e.~Theorem \ref{thm:FESp4} follows from the following important proposition in the same way that Theorem \ref{thm:FESL3} followed from Proposition \ref{prop:mainSL3}.
\begin{prop} \label{prop:xi_lambda_Cauchy_Sp4}
There is a vector $\xi_\infty$ such that for every cone $K\subset (0,\infty)^8$ with compact basis, 
\[
    \lim_{\lambda \in K, \, \lambda \to \infty} \|\xi_\lambda-\xi_\infty\|=0.
\]
\end{prop}
\begin{proof}[Proof of Theorem~\ref{thm:FESp4} assuming Proposition~\ref{prop:xi_lambda_Cauchy_Sp4}]
Since the vectors $\xi_{aaaaaaaa}$ are almost $X_1(\R)$-invariant, we deduce that the action $\Sp_4(\R)  \curvearrowright E$ is not proper, so it has a fixed point by \cite{MR2541757} and \cite[Lemma 2.14]{MR2316269}.
\end{proof}

Again, an important part of the proof will be of combinatorial nature, namely that a certain graph $\mathcal G$ is connected. For reasons of brevity, we first introduce this graph and then state the estimates needed in order to show that it is connected. Before defining $\mathcal{G}$, we first introduce three simpler graphs $\mathcal{G}^1$, $\mathcal{G}^2$, and $\mathcal{G}^0$.
\begin{defn} \label{defn:graphs_Sp4}
Let $\mathcal{G}^1=(V^1,A^1)$ be the graph with vertex set $V^1=(0,\infty)^4$ and an edge between $abcd$ and $ab'cd$ if $\max(b,b') < \min(a+d,\frac{a+c}{2})$ and between $abcd$ and $abc'd$ if $\max(c,c') < b+d$.

Let $\mathcal{G}^2=(V^2,A^2)$ be the graph with vertex set $V^2=(0,\infty)^4$ and an edge between $dcba$ and $dcb'a$ if and only if there is an edge in $\mathcal{G}^1$ between $abcd$ and $ab'cd$.

Let $\mathcal{G}^0=(V^0,A^0)$ be the graph with vertex set $V^0=(0,\infty)^8$ and the following edges:
\begin{enumerate}[(a)]
    \item\label{item:edge_of_type_odd} For $i=1,3,5$: $a_1a_2a_3a_4a_5a_6a_7a_8 \sim a_1'a_2'a_3'a_4'a_5'a_6'a_7'a_8'$ if $a_k'=a_k$ for all $k \leq i$ and $k \geq i+3$ and $a_ia_{i+1}a_{i+2}a_{i+3} \sim_1 a_i'a_{i+1}'a_{i+2}'a_{i+3}'$ (where $\sim_1$ denotes adjacency in $\mathcal{G}^1$).
    \item\label{item:edge_of_type_even} For $i=2,4$: $a_1a_2a_3a_4a_5a_6a_7a_8 \sim a_1'a_2'a_3'a_4'a_5'a_6'a_7'a_8'$ if and only if $a_k'=a_k$ for all $k \leq i$ and $k \geq i+3$ and $a_ia_{i+1}a_{i+2}a_{i+3} \sim_2 a_i'a_{i+1}'a_{i+2}'a_{i+3}'$ (where $\sim_2$ denotes adjacency in $\mathcal{G}^2$).
\end{enumerate}
Now, let $\mathcal{G}=(V,A)$ be the graph with vertex set $V=(0,\infty)^8 \times \{1,2\}$ and an edge between $(\lambda,i)$ and $(\lambda',j)$ if
\begin{enumerate}
    \item\label{item:noflipSp4_1} $i=j=1$ and $(\lambda,\lambda') \in A^0$,
    \item\label{item:noflipSp4_2} $i=j=2$ and $(\check\lambda,\check\lambda') \in A^0$,\\where $\check{\lambda} = a_8 a_7 a_6 a_5 a_4 a_3 a_2 a_1$ if $\lambda = a_1 a_2 a_3 a_4 a_5 a_6 a_7 a_8$,
    \item\label{item:flip_Sp4_condition} $i=1$, $j=2$ and $\lambda =  a_1 a_2 a_3 a_4 a_5 a_6 a_7 a_8$ and  $\lambda'=a_4 a_3 a_2 a_1 a_8 a_7 a_6 a_5$ with
    \[
        a_1+a_4<a_2, \; a_2+a_4<a_3, \; a_5+a_8<a_6, \; a_6+a_8<a_7.
    \]
\end{enumerate}
Finally, for $\varepsilon \in (0,1)$, denote by $\mathcal{G}^0_\varepsilon$ the graph obtained in the same way as $\mathcal{G}^0$ by replacing the strict inequalities defining $\mathcal{G}^1$ by the stronger inequalities $\max(b,b') \leq (1-\varepsilon) \min(a+d,\frac{a+c}{2})$ and $\max(c,c') \leq (1-\varepsilon)(b+d)$. Denote also $\mathcal{G}_\varepsilon$ the graph obtained in the same way as $\mathcal{G}$ by replacing the edges in $A_0$ by the edges in $A^0_\varepsilon$, and the conditions in \eqref{item:flip_Sp4_condition} by $a_1+a_4\leq(1-\varepsilon)a_2$, $a_2+a_4\leq(1-\varepsilon)a_3$, $a_5+a_8\leq(1-\varepsilon)a_6$, $a_6+a_8\leq(1-\varepsilon)a_7$.
\end{defn}

These definitions are motivated by the following consequence of Propositions~\ref{prop:nunutildechangeb} and \ref{prop:changec_H}. For $L>1$, we denote 
\[
    P_L = \cup_{s>0} [s,Ls]^8.
\]
\begin{prop}\label{prop:local_estimates_Sp4}
  For every $\varepsilon>0$ and $L>0$, there exists positive real numbers $C$ and $q(\varepsilon,L)<1$ such that for every $\lambda,\lambda' \in P_L$,
  \begin{equation}\label{eq:estimates_graph_1}
      (\lambda,1)\sim_\varepsilon (\lambda',1) \implies \| \xi_\lambda - \xi_{\lambda'}\| \leq C q(\varepsilon,L)^{\sqrt{|\lambda|}},
  \end{equation} 
  \begin{equation}\label{eq:estimates_graph_2}
 (\lambda,2)\sim_\varepsilon (\lambda',2) \implies \| \widetilde{\xi}_\lambda - \widetilde{\xi}_{\lambda'}\| \leq C q(\varepsilon,L)^{\sqrt{|\lambda|}},
   \end{equation} 
  \begin{equation}\label{eq:estimates_graph_12}
 (\lambda,1)\sim_\varepsilon (\lambda',2) \implies \| \xi_\lambda - \widetilde{\xi}_{\lambda'}\| \leq C q(\varepsilon,L)^{\sqrt{|\lambda|}}.  
 \end{equation} 
\end{prop}
\begin{proof}
The two inequalities \eqref{eq:estimates_graph_1} and \eqref{eq:estimates_graph_2} are proved in the same way, so let us focus on \eqref{eq:estimates_graph_1}. Suppose that $(\lambda,\lambda')$ is an edge in $A^0_\varepsilon$. Consider first an edge of type \eqref{item:edge_of_type_odd} for, say, $i=1$. Then $\lambda=a_1a_2a_3a_4a_5a_6a_7a_8$ and $\lambda'=a_1a'_2a'_3a_4a_5a_6a_7a_8$ with either $a'_2=a_2$ or $a'_3=a_3$. By the key observation, there is a continuous homomorphism $\varphi\colon H \to \Sp_4(\R)$ such that $\xi_\lambda = \varphi(\nu_{a_1a_2a_3a_4})\cdot \eta$ and $\xi_{\lambda'} = \varphi(\nu_{a_1a'_2a'_3a_4})\cdot \eta$ where $\eta = X_5(\gamma_{a_5})X_6(\gamma_{a_6})X_7(\gamma_{a_7})X_8(\gamma_{a_8}) \cdot \xi$. Observe that
\begin{align*}
    \max_{i \leq 4} \max_{|t| \leq e^{a_i}} & \|X_{i}(t) \cdot \eta - \eta\| \\ 
    &\leq \max_{i \leq 4} \max_{|t| \leq e^{|a_i|}}\|X_{i}(t) \cdot \xi- \xi\| &+ \sum_{i=5}^8 \left(1+\sqrt{\frac 2 \pi}\right) \max_{|t| \leq e^{a_i}}\|X_{i}(t) \cdot \xi- \xi\| \\ &\lesssim M(1+|\lambda|).
\end{align*}
In this inequality, we have used Lemma~\ref{lem:TV_of_gaussiansbis} and \eqref{eq:displacement_Sp4}. 

If $a'_3=a_3$, we therefore deduce from Proposition~\ref{prop:nunutildechangeb} that
\[
    \|\xi_\lambda - \xi_{\lambda'}\| \lesssim M(1+|\lambda|)q^{\sqrt{\varepsilon \min(a_1+a_4,\frac{1}{2}(a_1+a_3))}}.
\]
We deduce \eqref{eq:estimates_graph_1} from the fact that on $P_L$, we have $\min(a_1+a_4,\frac{1}{2}(a_1+a_3)) \geq \frac{1}{L}|\lambda|$.

If $a_2=a'_2$, we obtain the same inequality by applying Proposition~\ref{prop:changec_H}. For the other values of $i$ or an edge of type \eqref{item:edge_of_type_even}, the same argument applies.
    
Similarly, \eqref{eq:estimates_graph_12} follows by two applications of Lemma~\ref{lem:flip_H}.
\end{proof}

The following is the main combinatorial result in this section, the proof of which will be postponed to Section \ref{subsec:Sp4combinatorial}.
\begin{prop} \label{prop:connectednessSp4}
The graph $\mathcal{G}$ is connected.
\end{prop}

Proposition \ref{prop:xi_lambda_Cauchy_Sp4} now follows in the same way as Proposition \ref{prop:mainSL3}, using the graph $\mathcal{G}_{\varepsilon}$ as introduced in this section (instead of the corresponding graph for $\SL_3(\R)$). As for $\SL_3(\R)$, the full statement of Proposition~\ref{prop:connectednessSp4} is not needed for Theorem~\ref{thm:FESp4}, and it is enough to find $\lambda \in (0,\infty)^8$, $c>1$ and a path between $(\lambda,1)$ and $(c\lambda,1)$ in $\mathcal{G}$. An example of such a path of length $14$ with $c=1.125$ is given in Figure~\ref{figure:exampleSp4}.

\begin{figure}\label{figure:exampleSp4}
\begin{tikzpicture}
  \newcommand{\bighistogram}[9]{% #1-#8: heights, #9: position (e.g., (x,y))
    \begin{scope}[shift={#9}]
    \pgfmathsetmacro{\barwidth}{0.27} % Width of each bar
    \pgfmathsetmacro{\totalwidth}{8 * \barwidth} % Total width of the histogram
    \pgfmathsetmacro{\offset}{-\totalwidth/2 + \barwidth/2} % Center the histogram

    \foreach \i/\height in {1/#1, 2/#2, 3/#3, 4/#4, 5/#5, 6/#6, 7/#7, 8/#8} {
      \pgfmathsetmacro{\xpos}{\offset + (\i-1)*\barwidth} % Compute x-position explicitly
      \fill[blue!70] (\xpos, 0) rectangle ++(\barwidth, \height/30);
      % Add height label below the bar
      \node[below, font=\tiny] at (\xpos + \barwidth/2, -0.1) {\pgfmathprintnumber{\height}};
    }
  \end{scope}
}

  \pgfmathsetmacro{\deltax}{3}
  \pgfmathsetmacro{\deltay}{3}
  \pgfmathsetmacro{\xpos}{0}
  \pgfmathsetmacro{\ypos}{0}

  \newcommand{\deplacegauche}[1]{
    \pgfmathsetmacro{\xstart}{\xpos-1}
    \pgfmathsetmacro{\xend}{\xpos-\deltax+1.1}
    \draw[->] (\xstart,\ypos+0.5) -- node[above,midway] {#1} (\xend,\ypos+0.5);
    \pgfmathsetmacro{\xpos}{\xpos-\deltax}
    }

    \newcommand{\deplacedroite}[1]{
    \pgfmathsetmacro{\xstart}{\xpos+1.1}
    \pgfmathsetmacro{\xend}{\xpos+\deltax-1}
    \draw[->] (\xstart,\ypos+0.5) -- node[above,midway] {#1} (\xend,\ypos+0.5);
    \pgfmathsetmacro{\xpos}{\xpos+\deltax}
    }
    \newcommand{\deplacebas}[1]{
    \pgfmathsetmacro{\ystart}{\ypos-0.5}
    \pgfmathsetmacro{\yend}{\ypos-\deltay+1.6}
    \draw[->] (\xpos,\ystart) -- node[left,midway] {#1} (\xpos,\yend);
    \pgfmathsetmacro{\ypos}{\ypos-\deltay}
    }

  % Draw the 9 histograms in a row manually

\bighistogram{8}{24}{40}{8}{8}{24}{40}{8}{(\xpos,\ypos)}
\draw[dashed] (\xpos+\deltax/2, \deltay/2) -- ((\xpos+\deltax/2, -2.2*\deltay);
\deplacedroite{flip}
\bighistogram{8}{40}{24}{8}{8}{40}{24}{8}{(\xpos,\ypos)}
\deplacebas{4+5+}
\bighistogram{8}{40}{24}{9}{9}{40}{24}{8}{(\xpos,\ypos)}
\deplacegauche{flip}
\bighistogram{9}{24}{40}{8}{8}{24}{40}{9}{(\xpos,\ypos)}
\deplacegauche{4+3+}
\bighistogram{9}{24}{46}{23}{8}{24}{40}{9}{(\xpos,\ypos)}
\deplacegauche{2+5+}
\bighistogram{9}{27}{46}{23}{35}{24}{40}{9}{(\xpos,\ypos)}
\deplacebas{6+7+}
\bighistogram{9}{27}{46}{23}{35}{37}{45}{9}{(\xpos,\ypos)}
\deplacedroite{3-6-}
\bighistogram{9}{27}{45}{23}{35}{27}{45}{9}{(\xpos,\ypos)}
\deplacedroite{5-4-}
\bighistogram{9}{27}{45}{9}{9}{27}{45}{9}{(\xpos,\ypos)}

\end{tikzpicture}

    \centering
\caption{A path of length $14$ between $(8\lambda,1)$ and $(9\lambda,1)$ in $\mathcal{G}$ for $\lambda=(1,3,5,1,1,3,5,1)$. The dashed line represents the boundary between $V^0\times \{1\}$ on the left and $V^0 \times \{2\}$ on the right. A vertex $(\lambda,i) \in V$ is represented by a histogram with heights $\lambda$. An edge coming from $\mathcal{G}^0$ is oriented and labeled by the corresponding position and a sign, depending on whether the coordinate is increased or decreased. The other edges are labeled \emph{flip}. To shorten the picture, an edge with two labels means the concatenation of two edges. For example, the vertical edge labeled $4+5+$ means the path $((8, 40, 28,8,8,40,24,8),2) - ((8, 40, 28,9,8,40,24,8),2)  - ((8, 40, 28,9,9,40,24,8),2)$.}
\end{figure}

\subsection{Combinatorial part} \label{subsec:Sp4combinatorial}
We are left to prove Proposition~\ref{prop:connectednessSp4}. The graphs $\mathcal G,\mathcal G^0,\mathcal G^1,\mathcal G^2$ have been defined in Definition~\ref{defn:graphs_Sp4}. The main ingredient is the following lemma.
\begin{lemma} \label{lemma:connectednessG0}
  Two vertices $\lambda=a_1a_2a_3a_4a_5a_6a_7a_8$ and $\lambda'=a_1'a_2'a_3'a_4'a_5'a_6'a_7'a_8'$ are in the same connected component of $\mathcal{G}^0$ if and only if $a_1=a_1'$ and $a_8=a_8'$.
\end{lemma}
We first establish the following sublemma.
\begin{lemma} \label{lem:zigzaggraph}
    If $abcd$ and $ab'c'd$ satisfy $\max(b,b') < a+d$ and $\max(c,c') < a+2d$, then $abcd$ and $ab'c'd$ are in the same component of $\mathcal{G}^1$. In particular, if $\max(b,c) \leq \max(a,d) = m$, then $abcd \sim_1 ammd$.
\end{lemma}
\begin{proof}
    Fix $a$ and $d$, and consider the induced subgraph of $\mathcal{G}_1$ with vertex set $\{abcd\mid b<a+d, c<a+2d\}$. Making the change of variable $x=a+d-b$ and $y=a+2d-c$, this induced subgraph is isomorphic to the graph on $(0,a+d)\times(0,a+2d)$ with edges
    \begin{align*} xy \sim x'y &\iff \min(x,x') > \frac y 2,\\
    xy \sim xy' & \iff \min(y,y')> x.\end{align*} 
    It is clear (see Figure~\ref{picture:induced_graph}) that this subgraph is connected. This proves the lemma.
    \begin{figure}
  \center
\begin{tikzpicture}
\begin{scope}[scale=2]
\draw  (0,0)--(2,0)--(2,3)--(0,3)--(0,0);
\draw[dotted] (0,0)--(2,2);
\draw[dotted] (0,0)--(1.5,3);
\foreach \x in {0,...,6}
 {\draw[blue] (0.25*\x,0.5*\x)--(2,0.5*\x);
 \draw[darkgreen] (0.3*\x,0.3*\x)--(0.3*\x,3);}
\end{scope}
\end{tikzpicture}
\caption{The induced subgraph of $\mathcal{G}^0$ from Lemma~\ref{lem:zigzaggraph}.}\label{picture:induced_graph}
\end{figure}

\end{proof}

\begin{proof}[Proof of Lemma \ref{lemma:connectednessG0}]
It is clear from the definition that if $a_1a_2a_3a_4a_5a_6a_7a_8$ is in the same connected component as $a'_1a'_2a'_3a'_4a'_5a'_6a'_7a'_8$, then $a_1=a'_1$ and $a_8=a'_8$.

For the converse, fix $a_1,a_8$. From Lemma~\ref{lem:zigzaggraph}, it follows that $a_1a_2a_3a_4a_5a_6a_7a_8$ is connected to $a_1 mmmmmm a_8$ for $m=\max a_i$. Now if $m \geq a_1,a_8$ and $m<m'<m+a_1$, by Lemma~\ref{lem:zigzaggraph}, we know that $a_1mmmmmm a_8$ is connected to $a_1 m m' mmmm a_8$, which is connected to $a_1m'm'm'm'm'm' a_8$ by the previous discussion. As a consequence, the vertices $\{ a_1 mmmmmm a_8 \mid m \geq \max(a_1,a_8)\}$ are all in the same connected component. Moreover, every element of the form $a_1a_2a_3a_4a_5a_6a_7a_8$ is connected to one of these vertices. This proves the lemma.
\end{proof}
\begin{proof}[Proof of Proposition~\ref{prop:connectednessSp4}]
It suffices to prove that for every $\lambda,\lambda' \in (0,\infty)^8$, the vertex $(\lambda,1)$ is in the same connected component as $(\lambda',2)$. Pick $x>\max( a_1,a'_1,a_8,a'_8)$. Then by Lemma~\ref{lemma:connectednessG0}, the vertex $(\lambda,1)$ is in the same connected component as $(a_1\ 2x\ 3x\ a'_1\ a'_8\ 2x\ 3x\ a_8,1)$, and $(\lambda',2)$ is in the same connected component as $(a'_1\ 3x\ 2x\ a_1\ a_8\ 3x\ 2x\ a'_8,2)$. Our choice of $x$ guarantees that these two vertices are connected by an edge of type \eqref{item:flip_Sp4_condition} in $\mathcal{G}$.
\end{proof}

\subsection{Extension to the universal cover of $\Sp_4(\R)$}
In this short section, we point out that a large part of Theorem~\ref{thm:FESp4} remains true for the universal cover $\widetilde{Sp}_4(\R)$ of $\Sp_4(\R)$. Note that this is an infinitely sheeted covering. In a forthcoming work, we will apply this to extend Theorem~\ref{thm:higherrankgroupsFE} to infinite center groups.
\begin{thm}\label{eq:FE_Sp4tilde}
Every isometric action of $\widetilde{\Sp}_4(\R)$ on a uniformly convex Banach space has a point fixed under the action of $\exp(\R \mathfrak{X}_{2\beta})$.
\end{thm}
Indeed, as above, for every root $\phi$ in the root system $C_2$, the exponential map $\exp \colon \mathfrak{sp}_4 \to \widetilde{Sp}_4(\R)$ gives rise to a continuous homomorphism from $\R$ to $\widetilde{Sp}_4(\R)$ given by $X_\phi(t) = \exp(t \mathcal{X}_\phi)$. Using that the group $H$ is simply connected, we see that any four consecutive roots in the root system $C_2$ give rise to an embedding of $H$ into $\widetilde{Sp}_4(\R)$. For $\lambda \in (0,\infty)^8$, we can define measures $\mu_\lambda$ and $\widetilde{\mu}_\lambda$ on $\widetilde{Sp}_4(\R)$ by exactly the same formulas as for $\Sp_4(\R)$.

If $G$ acts by isometries on a uniformly convex Banach space $E$ and $\xi$ is an element of $E$, we define, in the same way, $\xi_\lambda = \mu_\lambda \cdot \xi$ and $\widetilde \xi_\lambda = \widetilde{\mu}_\lambda \cdot \xi$. Lemma~\ref{lem:orbitgrowth} implies that there exists $M>0$ such that \eqref{eq:displacement_Sp4} holds. Therefore, every step of the proof of Proposition~\ref{prop:xi_lambda_Cauchy_Sp4} works identically in this setting, so Proposition~\ref{prop:xi_lambda_Cauchy_Sp4} is also true in this generality, and the vector $\xi_\infty$ is easily seen to be $X_1(\R)$-invariant.
\section{Proofs of Theorems \ref{thm:higherrankgroupsFE} and \ref{thm:BFGMconjecture}} \label{sec:conclusion}

\begin{proof}[Proof of Theorem \ref{thm:higherrankgroupsFE}]
Let $G$ be a connected simple Lie group with real rank $\geq 2$ and finite center, and consider an isometric action $G \curvearrowright E$ on a uniformly convex Banach space $E$. It is well known that $G$ has a closed subgroup $G'$ with finite center such that $G'/Z(G')$ is isomorphic to $\SL(3,\R)$ or to $\Sp(4,\R) / \{I_4,I_4\}$; see e.g.~\cite[Theorem I.1.6.2]{MR1090825}. Since property (F$_E$) passes to quotients (see \cite[Proposition 2.15.(2)]{MR2316269}) and extensions by finite groups (see \cite[Proposition 2.5.4]{MR2415834}, the proof of which carries over to property (F$_E$)), it follows from Theorem \ref{thm:FESL3} and Theorem \ref{thm:FESp4} that $G'$ has (F$_E$). It follows that the action $G \curvearrowright E$ cannot be proper, so (in a similar fashion as in Theorem \ref{thm:FESL3} and Theorem \ref{thm:FESp4}) we deduce that $G \curvearrowright E$ has a fixed point by \cite[Theorem 1.4]{MR2541757} and \cite[Lemma 2.14]{MR2316269}.
%It is well-known that its Lie algebra contains a subalgebra isomorphic to $\mathfrak{sl}_3(\R)$ or to $\mathfrak{sp}_4(\R)$, see e.g.~\cite[Theorem I.1.6.2]{MR1090825}. So it follows from Remark~\ref{rem:FE_containingSL3} or Theorem~\ref{eq:FE_Sp4tilde} that $G$ has F$_E$ for every uniformly convex Banach space $E$.

The statement for lattices in $G$ follows directly from \cite[Proposition 8.8]{MR2316269}, since every lattice in $G$ is $p$-integrable \cite{MR1767270} and the class of uniformly convex Banach spaces is stable under $E \mapsto L^p(G/\Gamma;E)$ for every $1 < p < \infty$.
\end{proof}

\begin{proof}[Proof of Theorem \ref{thm:BFGMconjecture}]
Let $G = \prod_{i=1}^n G_i (\F_i)$ be a higher rank group, i.e.~each $\F_i$ is a local field and each $G_i$ is a (Zariski) connected (almost) $\F_i$-simple group with $\F_i$-rank $\geq 2$. We need to show that $G$ and its lattices have property F$_E$ for every super-reflexive Banach space $E$. By \cite[Proposition 2.13]{MR2316269}, we can reduce to the case of isometric actions on uniformly convex spaces.

We first consider actions of $G$. Note that by \cite[Proposition 2.15.(3)]{MR2316269}, it suffices to prove that each factor $G_i(\F_i)$ has (F$_E$) for every uniformly convex Banach space. If $\F_i$ is Archimedean (i.e.~$\R$ or $\C$), this is Theorem \ref{thm:higherrankgroupsFE}. If $\F_i$ is non-Archimedean, this follows from the fact that $G_i(\F_i)$ has Lafforgue's strong Banach property (T); see \cite{MR2423763, MR2574023, MR3190138}.

Now, let $\Gamma$ be a lattice in $G$. It is well known (see e.g.~\cite[Theorem V.5.22]{MR0507234}) that $\Gamma$ is isomorphic (up to finite index) to a direct product of irreducible lattices in $G$. The statement for lattices now follows (in the same way as in the proof of Theorem \ref{thm:higherrankgroupsFE}) from \cite[Proposition 8.8]{MR2316269}, noting that every irreducible lattice in $G$ is $p$-integrable for every $1 < p < \infty$ \cite{MR1767270}.
\end{proof}

\begin{rem}
Note that the conjecture in \cite{MR2316269} is stated in terms of property $\overline{\mathrm{F}}_E$ rather than property F$_E$. Property $\overline{\mathrm{F}}_E$ is formally stronger than property F$_E$, but by \cite[Proposition 2.13]{MR2316269}, the statements of Theorem A and \cite[Conjecture 1.6]{MR2316269} are equivalent.
\end{rem}

\section{Spectral gap and super-expanders} \label{sec:superexpanders}

Expanders are sequences of finite, highly connected, sparse graphs with an increasing number of vertices. The notion of super-expander, which we recall now, was introduced by Mendel and Naor \cite{MR3210176} and refers to sequences of graphs that satisfy some nonlinear form of spectral gap with respect to all super-reflexive Banach spaces.
\begin{defn}
A sequence of finite, $d$-regular graphs $\mathcal{G}_n=(V_n,A_n)$, $n \in \mathbb{N}$, with $\lim_{n \to \infty} |V_n|=\infty$, is called an expander with respect to the Banach space $E$ if there exists $\gamma > 0$ such that for all $n \in \mathbb{N}$ and all $f:V_n \to E$, the following inequality holds:
\[
	\frac{1}{|V_n|^2}\sum_{v,w \in V_n} \| f(v)-f(w) \|_E^2 \leq \frac{\gamma}{d|V_n|} \sum_{v \sim w} \| f(v)-f(w) \|_E^2.
\]
A sequence of finite, $d$-regular graphs is called a super-expander if the above is true for all super-reflexive Banach spaces $E$.
\end{defn}
The usual notion of expander corresponds to taking $E=\C$, or equivalently, $E=\ell^2$. 
\begin{defn}
An isometric representation $\pi$ of a locally compact group $G$ on a Banach space $E$ is said to have spectral gap if there exists a compact subset $Q \subset G$ and a constant $c>0$ such that for all $\xi \in E$,
\[
    \max_{g \in Q} \|\pi(g) \xi- \xi\| \geq c \, d(\xi, E^\pi),
\]
where $d(\xi,E^\pi)$ denotes the distance between $\xi$ and the subspace $E^\pi = \{\eta \in E \mid \forall g \in G: \; \pi(g) \eta=\eta\}$ of $\pi(G)$-invariant vectors.
\end{defn}
If $G$ admits a compact generating set $Q_0$, then we can always take $Q=Q_0$ (because by Baire's theorem every compact subset of $G$ is contained in $(\{1\}\cup Q_0)^N$ for some $N$).
\begin{defn}
A locally compact group has property (T$_E$) if every isometric representation on $E$ has spectral gap.
\end{defn}

As mentioned in Section \ref{subsec:strategy}, by \cite[Theorem 1.3]{MR2316269}, we obtain the following consequence of Theorem~\ref{thm:BFGMconjecture}.
\begin{cor}\label{cor:TE}
    Higher rank groups and their lattices have property (T$_E$) for every super-reflexive Banach space $E$.
\end{cor}

Particularly interesting representations on super-reflexive Banach spaces arise as Koopman representations on vector-valued Bochner spaces of actions by measure-preserving transformations. More generally, if $G$ acts continuously on a standard measure space $(X,\mu)$ by measure-class preserving transformations, we obtain a continuous isometric representation of $G$ on $L^2(X,\mu;E)$ by the formula
\[
    \pi(g)f(x) = \left(\frac{ d g_*\mu}{d\mu}(x)\right)^{\frac 1 2} f(g^{-1}x).
\]
Since $L^2(X,\mu;E)$ is super-reflexive whenever $E$ is super-reflexive, a particular case of Corollary~\ref{cor:TE} is that the representation on $L^2(X,\mu;E)$ has spectral gap. For the special case of lattices and permutation representations, we obtain the following consequence, which contains Corollary~\ref{cor:superexpanders} as the particular case of finite index normal subgroups.
\begin{cor}
Let $\Gamma$ be a higher rank lattice and $S \subset \Gamma$ a finite generating set. For every super-reflexive Banach space $E$, there exists $C>0$ such that for every subgroup $\Lambda <\Gamma$ and every $f \in \ell^2(\Gamma/\Lambda;E)$,
\[
    \sum_{x,y \in \Gamma/\Lambda, \, x \in S y} \|f(x) - f(y)\|^2 \geq C \inf_{\eta \in E} \sum_{x \in \Gamma/\Lambda} \| f(x) -\eta\|^2.
\]
\end{cor}

\end{document}